\documentclass{amsart}
\usepackage{amsthm}
\usepackage{amsfonts}
\usepackage{amssymb}
\usepackage{amsmath}
\usepackage{mathrsfs}
\usepackage{csquotes}

\begin{document}
	
\large
	
\theoremstyle{plain}
\theoremstyle{definition}
\newtheorem{theorem}{Theorem}[section]
\newtheorem{prop}[theorem]{Proposition}
\newtheorem{lemma}[theorem]{Lemma}
\newtheorem{Corollary}{Corollary}[theorem]

\newtheorem{example}[theorem]{Example}

\newtheorem{definition}[theorem]{Definition}
\newtheorem{remark}[theorem]{\textnormal{\textbf{Remark}}}

\numberwithin{equation}{section}




%






\title[Roughness of nets]{Roughness of nets}

 \author[Sanjoy Ghosal and Mandobi Banerjee]{Sanjoy Ghosal$^{*}$ and Mandobi Banerjee$^{**}$}
 
 \thanks{$^{*}$Department of Mathematics, University of North Bengal, Raja Rammohunpur, Darjeeling-734013, West Bengal, India. \emph{E-mail}: sanjoykumarghosal@nbu.ac.in; sanjoyghosalju@gmail.com. ORCID: 0000-0001-8563-1941 (S. Ghosal).}

  \thanks{$^{**}$Department of Engineering Science and Humanities (Mathematics), Siliguri Institute of Technology, Hill Cart road, Sukna, Siliguri, Darjeeling-734009, West Bengal, India.	\emph{E-mail}: mandobibanerjee.desh@sittechno.edu.in; dr.banerjeemandobi@gmail.com. ORCID: 0000-0002-9790-2515 (M. Banerjee).}

\date{}














\date{}

\setcounter {page}{1}

\maketitle




\begin{abstract}  In this paper, we introduce the concepts of ~`\emph{rough weighted ideal-$\tau$-convergence}' and ~`\emph{rough weighted ideal-$\tau$-cluster point}' of a net in locally solid Riesz spaces. We also compare these concepts with key findings from the works of Lahiri $\&$ Das (2007/2008), Das (2012) and Aytar (2008), as follows:
		\begin{itemize}				
				\item[(i):] The set of \emph{ideal-cluster points} of a sequence in a topological space is closed, 
				\item[(ii):] The set of \emph{ideal-cluster points} of a net in a topological space is closed, 
				\item[(iii):] The set of \emph{rough ideal-cluster points} of a sequence in a norm linear space is closed. 
			\end{itemize}			
			We demonstrate that, generally, the set of ~`\emph{rough weighted ideal-$\tau$-cluster points}' of a net in a locally solid Riesz space may not satisfy the aforementioned results, which implies that it may not be closed. This suggests that this set exhibits new characteristics, which we will discuss here. Ultimately, our characterization demonstrates that these new concepts can help identify an ideal that is maximal and effectively represent any closed $\tau$-neighborhood of $\theta$ within these spaces.
\end{abstract}

\section{Introduction}

In the year 2008, Aytar \cite{A1, A2} introduced the notion of \emph{rough statistical convergence} as a generalization of \emph{rough convergence} \cite{Phu1,Phu2} and \emph{statistical convergence} \cite{Salat,Steinhaus}. He demonstrated \emph{rough statistical convergence} in the following manner: A sequence $x=\{x_{n}\}_{n\in \mathbb{N}}$ in a norm linear space $X$ is said to be \emph{rough statistically convergent} to $x_*$ w.r.t the roughness of degree $r$ (where $r\geq 0$) if for every $\varepsilon>0,$
\begin{equation}\label{Eq1.1}
\displaystyle{\lim_{n\rightarrow \infty}}~\frac{1}{n}|\{k\leq {n}: ||x_k-x_*||\geq r+\varepsilon\}|=0,
 \end{equation}   
and we write $x_{n}\xrightarrow[r]{st} x_*.$ Thus, the set $st-\mbox{LIM}^rx=\left\{x_*\in X: x_{n}\xrightarrow[r]{st} x_* \right\}$ is called the \emph{rough  statistical limit set} of the sequence $x=\{x_{n}\}_{n\in \mathbb{N}}$ with degree of roughness $r.$ Moreover Aytar \cite{A3} introduce the concept of \emph{rough statistical cluster point} in the following way: Let $r$ be a non-negative real number. An element $\zeta$ in a norm linear space is called a  \emph{rough statistical cluster point} w.r.t the roughness degree $r$ (or, \emph{r-statistical cluster point}) of a sequence $x=\{x_{n}\}_{n\in \mathbb{N}}$  provided that
 \begin{equation}\label{Eq1.2} \displaystyle{\lim_{n\rightarrow \infty}}~\frac{1}{n}|\{k\leq {n}: ||x_k-\zeta||< r+\varepsilon\}|\neq 0, \end{equation} 
for every $\varepsilon>0.$  The set of all \emph{r-statistical cluster points} of the sequence $x=\{x_{n}\}_{n\in \mathbb{N}}$ is symbolized by $\Gamma^{r}_x.$ If $r=0$ in the above expression we could visualize \emph{r-statistical cluster point} exactly to be merged with ordinary statistical cluster point \cite{Fridy, fridy1993statistical}.

\begin{theorem} \label{fridy} \cite[Proposition 2]{fridy1993statistical} For any sequence of real numbers $x=\{x_n\}_{n\in \mathbb{N}},$ the set $\Gamma_x$ of statistical cluster points of the sequence $x=\{x_n\}_{n\in \mathbb{N}}$ is a closed set. \end{theorem}

\begin{theorem} \label{Di Maio et al.} \cite[Theorem 2.9]{Di} In a topological space the statistical cluster points set of a sequence $x=\{x_n\}_{n\in \mathbb{N}}$ is closed. \end{theorem}

\begin{theorem} \label{1} \cite[Theorem 1.1]{A3} For any sequence $x=\{x_n\}_{n\in \mathbb{N}}$ in a norm linear space $X,$ the rough statistical cluster point set $\Gamma^{r}_x$ is closed, for every $r\geq 0.$\end{theorem}

Di Maio and Ko$\check{\mbox{c}}$inac et al. \cite{Di} introduced the concept of statistical convergence in topology. Our aim is to introduce rough weighted ideal convergence and rough weighted ideal cluster points for a net in the context of topological vector spaces. In this context, we introduce the following definitions and results in sequence, which are primarily required.\\

Dündar et al. \cite{Dundar2014,Dundar2016} and Pal et al. \cite{pal} independently introduced the concept of \emph{rough $\mathcal{I}$-convergence} as a generalization of \emph{rough convergence} \cite{Phu1,Phu2}, \emph{rough statistical convergence} \cite{A1,LR1,Garcia}, and \emph{$\mathcal{I}$-convergence} \cite{Das-Savas2014(Taiwan),Kos1, Kos2,SGSG(Taiwan)}. This notion is based on the structure of the \emph{ideal} $\mathcal{I}$ of subsets of the set $\mathbb{N}$. A family $\mathcal{I}\subset 2^{\mathbb{N}}$ is called an \emph{ideal} on $\mathbb{N}$ if (i) $\varnothing \in \mathcal{I}$ (ii) $A,B\in \mathcal{I}$ implies $A\cup B\in \mathcal{I}$ (iii) $A\in \mathcal{I}, B\subset A$ imply $B\in \mathcal{I}.$ The ideal $\mathcal{I}$ is called \emph{non-trivial} if $\mathcal{I}\neq \{\varnothing\}$ and $\mathbb{N}\notin \mathcal{I}$ and a non-trivial ideal $\mathcal{I}$  is called \emph{admissible} if $\mathcal{I}$ contains all the singleton set.  Clearly $\mathcal{I}\subset 2^{\mathbb{N}}$ is a non-trivial ideal of $\mathbb{N}$ iff $\mathcal{F}=\mathcal{F}(\mathcal{I})=\{K\subset \mathbb{N}:\mathbb{N}\setminus K\in \mathcal{I}\}$ is a \emph{filter} on $\mathbb{N},$ called the \emph{filter associated with $\mathcal{I}.$} For more details on the various types of rough convergence in recent work, see references \cite{Aziz-Ghosal,Leonetti, Mandal-Ghosal,mursaleen}.

\begin{definition} \cite{Dundar2014,Dundar2016,pal} Let $r$ be a non-negative real number. A sequence  $x=\{x_n\}_{n\in \mathbb{N}}$ in a norm linear space is said to be \emph{rough $\mathcal{I}$-convergent} to $x_*$ w.r.t the roughness of degree $r$ (where $r\geq 0$), denoted by $x_n\xrightarrow[r]{\mathcal{I}}x_*,$ provided for any $\varepsilon>0,$ the set $\{n\in \mathbb{N}: ||x_n-x_*||\geq r+\varepsilon\}\in \mathcal{I}.$ 
\end{definition}

\begin{definition} \cite{CSST, Kos2} Let $(X,d)$ be a metric space. A point $ \zeta \in X$ is called an \emph{$\mathcal{I}$-cluster point} of a sequence $x=\{x_n\}_{n\in \mathbb{N}}$ in $X$ if for any $\varepsilon>0,$ the set $\{n\in \mathbb{N}: d(x_n,\zeta)<\varepsilon\}\notin \mathcal{I}.$
\end{definition}
		
	The set of all \emph{$\mathcal{I}$-cluster point} of $x=\{x_n\}_{n\in \mathbb{N}}$ in a metric space $(X,d)$  is denoted by $\Gamma_x(\mathcal{I})$ and the set $\Gamma_x(\mathcal{I})$ is closed in $X,$ for any admissible ideal $\mathcal{I}$ \cite{BL,CSST, Kos2}.

\begin{theorem} \label{theorem 1.2} \cite{CSST} Let $\mathcal{I}$ be an admissible ideal on $\mathbb{N},$ $\Gamma_{x}(\mathcal{I})$ be the set of all $\mathcal{I}$-cluster points of a sequence $x=\{x_n\}_{n\in \mathbb{N}}$  and $K$ be a compact subset of a metric space $X.$ Then for every sequence $x=\{x_n\}_{n\in \mathbb{N}}$ in $X$ with $\{n\in \mathbb{N}: x_n\in K\}  \notin \mathcal{I}$ we have $K \cap \Gamma_{x}(\mathcal{I})\neq \varnothing.$ \end{theorem}

Lahiri et al. \cite{Lahiri} initiate the idea of \emph{$\mathcal{I}$-convergence net} in a topological space and establish several results using the concept. They define the \emph{$\mathcal{I}$-convergence net} as follows: Let $(D,\geq)$ be a directed set and $\alpha_0 \in D, ~D_{\alpha_0}=\{\alpha\in D: \alpha \geq \alpha_0\}.$ An ideal $\mathcal{I} \subset 2^{D}$ is called \emph{non-trivial} if $\mathcal{I}\neq \{\varnothing\}$ and $D\notin \mathcal{I}.$  A non-trivial ideal $\mathcal{I}$ of $D$ will be called \emph{$D$-admissible}, if $D_{\alpha}\in \mathcal{F(\mathcal{I})}$ for all $\alpha \in D.$ If $\mathcal{I}$ is a non-trivial ideal of $D,$ then $\mathcal{F(\mathcal{I})}=\{M\subset D: M^c\in \mathcal{I}\}$ is a \emph{filter associated with the ideal} $\mathcal{I}.$ For example the collection 
$$\mathcal{F}_0=\{A\subset D: ~\mbox{there exists a}~D_{\alpha_0}~\mbox{with}~A\supset D_{\alpha_0}\}$$ 
forms a filter in $D.$ Therefore $\mathcal{I}_0=\{A\subset D: A^c\in \mathcal{F}_0\}$ is a non-trivial \emph{$D$-admissible} ideal of $D.$

\begin{definition} \cite{Lahiri}  Let $\mathcal{I}$ be a non-trivial ideal of the directed set $D.$ A net $s=\{s_{\alpha}\}_{\alpha\in D}$ in an arbitrary topological space $X$ is said to be \emph{$\mathcal{I}$-convergent} to $x_0\in X,$ if for any open set $\textbf{U}_{x_0}$ containing $x_0$ the set $\{\alpha \in D: s_\alpha\notin \textbf{U}_{x_0}\}\in \mathcal{I}.$ Symbolically we write $\mathcal{I}-\lim s_\alpha=x_0$ and we say that $x_0$ is the $\mathcal{I}$-limit of the net $s=\{s_{\alpha}\}_{\alpha\in D}.$
\end{definition}

\begin{definition} \cite{Lahiri}  Let $\mathcal{I}$ be a non-trivial ideal of the directed set $D$ and $s=\{s_{\alpha}\}_{\alpha\in D}$ be a net in an arbitrary topological space $X.$ A point $ \zeta \in X$ is called \emph{$\mathcal{I}$-cluster point} of the net $s=\{s_{\alpha}\}_{\alpha\in D},$ if for any open set $\textbf{U}_{\zeta}$ containing $\zeta$ the set $\{\alpha \in D: s_\alpha\in \textbf{U}_\zeta\}\notin \mathcal{I}.$ 
\end{definition}

\begin{theorem} \label{theorem 1.1} \cite[Theorem 9]{Lahiri} Let $\mathcal{I}$ be an ideal on a directed set $D.$ In a topological space, the {$\mathcal{I}$-cluster points set of a net} is closed. \end{theorem}

\begin{definition} \label{subnet} \cite{Lahiri} A net $\omega =\{\omega_{e}\}_{e\in E},$ where $E$ is a directed set, is said to be a subnet of the net $s=\{s_{\alpha}\}_{\alpha\in D}$ if there is a mapping $i: E \rightarrow D$ such that 	
	\begin{description}
		\item[(i)] $\omega=soi$
		\item[(ii)] for any $\alpha \in D$ there is an element $e_0\in E$ such that $i(e)\geq \alpha$ for all $e\in E$ with $e\geq e_0.$
	\end{description}	
\end{definition}

\begin{remark} \label{new ideal}
	Let $\omega =\{\omega_{e}\}_{e\in E}$ be a subnet of the net $s=\{s_{\alpha}\}_{\alpha\in D}$ and $\mathcal{I}^E$ be an $E$-admissible ideal of $E.$ Then the collection $\mathcal{I}^E(i)=\{A\subset D: i^{-1}(A)\in \mathcal{I}_E\}$ is a $D$-admissible ideal of $D.$ Since each $\alpha \in D$ there is an element $e_0\in E$ such that $i(e)\geq \alpha$ for all $e\in E$ with $e\geq e_0.$ This implies $i(E_{e_0})\subset D_{\alpha}$ where $E_{e_0}=\{e\in E: e\geq e_0\}.$ Thus $E_{e_0}\subset i^{-1}(D_{\alpha})$ and $E_{e_0} \in \mathcal{F}(\mathcal{I}^E).$ So $i^{-1}(D_{\alpha}) \in \mathcal{F}(\mathcal{I}^E).$ Consequently $D_{\alpha} \in \mathcal{F}(\mathcal{I}^E(i)).$
\end{remark}

 On the other hand, we briefly recall some of the basic definitions in the theory of \emph{Riesz space}, and we refer the readers to \cite{AP,R,SP} for more details.

\begin{definition} \cite{AP, GM2} Consider $L$ is a vector space over the field $\mathbb{R}$ as well as  a lattice with respect to the partial ordering $\leq$ on $L.$ Then $L$ is said to be a \emph{Riesz space} (or a \emph{vector lattice}) if it satisfies the following properties:\\
	(i) If $a, b \in L$ and $b \leq a,$ then $b + c \leq a+ c$ for each $c\in L,$\\
	(ii) If $a, b\in L$ and $b \leq a,$ then $\lambda b \leq \lambda a$ for each $\lambda \geq 0.$
\end{definition}

The \emph{absolute value} of $a\in L$ is defined by $|a| = a \vee(-a)$ and  $a^{+} = a\vee \theta,$ $a^{-} = (-a)\vee\theta,$ where $\theta$ is the  \emph{zero element} of $L.$ A subset $A$ of a \emph{Riesz space} $L$ is said to be \emph{solid}, if $b\in A$ and $|a| \leq |b|$ imply $a\in A.$ A topology $\tau$ on $L$ is a \emph{linear topology} if the addition and the scalar multiplication are continuous. Every \emph{linear topology} $\tau$ on a vector space $L$ has a base $\mathcal{N}$ for the $\tau$-neighbourhoods of $\theta$ (zero) satisfying the following properties:

\begin{itemize}
	\item[(i):] Each $\textbf{U}_\theta\in\mathcal{N}$ is a \emph{balanced set},
	
	\item[(ii):] Each $\textbf{U}_\theta \in\mathcal{N}$ is an \emph{absorbing set},
	
    \item[(iii):] For each $\textbf{U}_\theta\in \mathcal{N}$ there exists some $\textbf{W}_\theta \in \mathcal{N}$ with $\textbf{W}_\theta + \textbf{W}_\theta \subseteq \textbf{U}_\theta.$   
\end{itemize}

\begin{definition} \cite{AP, GM2} A \emph{linear topology} $\tau$ on a \emph{Riesz space} $L$ is said to be \emph{locally solid} if $\tau$ has a base at $\theta$ consisting of \emph{solid sets}. A \emph{locally solid Riesz space} $(L, \tau)$ is a \emph{Riesz space} $L$ equipped with a \emph{locally solid topology} $\tau.$ \end{definition}

The symbol $\mathcal{N}_{sol}$ stands for a base at $\theta$ consisting of \emph{solid sets} and satisfying the properties (i), (ii), and (iii) in a \emph{locally solid topology}.

\begin{definition} \cite{GB, GM2} A subset $A$ of a \emph{Riesz space} $(L,\tau)$ is said to be \emph{$\tau$-bounded} (or \emph{topologically bounded}) if for every $\tau$-neighbourhood $\textbf{U}_\theta$ of $\theta$ there exists some $\lambda>0$ such that $A\subset \lambda \textbf{U}_\theta.$ If $A$ is not \emph{$\tau$-bounded}, then it is called \emph{$\tau$-unbounded}.
\end{definition}

\begin{definition} \cite{DS(filomat),GM2} Let $(L,\tau)$ be a \emph{locally solid Riesz space}. A net $s=\{s_{\alpha}\}_{\alpha\in D}$ in $L$ is said to be \emph{$\mathcal{I}_{\tau}$-convergent} (or, \emph{ideal-$\tau$-convergent}) to $x_0\in L$ if for any $\tau$-neighbourhood $\textbf{U}_\theta$ of $\theta,$ $\{\alpha \in D: s_\alpha-x_0\notin \textbf{U}_\theta\}\in \mathcal{I}.$ Briefly we write $\mathcal{I}_\tau-\lim s_\alpha=x_0$ (or, $s_\alpha\xrightarrow{\mathcal{I}_\tau} x_0$).
\end{definition}


\begin{definition} \cite{DS(filomat),GM2} A net $s=\{s_{\alpha}\}_{\alpha\in D}$ in a locally solid Riesz space $(L,\tau)$ is said to be \emph{$\mathcal{I}_{\tau}$-bounded} (or, \emph{ideal-$\tau$-bounded}), if for any $\tau$-neighbourhood $\textbf{U}_\theta$ of $\theta$ there exists some $\lambda>0$ such that $\{\alpha\in D: \lambda s_{\alpha}\notin \textbf{U}_\theta\}\in \mathcal{I}.$   \end{definition}

In this discourse, we now proceed to the definitions of \emph{rough weighted $\mathcal{I}_{\tau}$-convergence} \cite{GM1} and \emph{weighted $\mathcal{I}_{\tau}$-cluster point} \cite{GM1} of a sequence in a \emph{locally solid Riesz space} $(L,\tau).$ Meanwhile, we define the \emph{weighted sequence} and  \emph{weighted convergence sequence} of real numbers.


\begin{definition} \label{weight sequence} \cite{GB,GM1,GM2} A  sequence of real numbers $t=\{t_{n}\}_{n\in \mathbb{N}}$ is said to be weighted sequence if there exists a positive real number $\delta$ such that
	\begin{equation} \label{Eq(1)} t_n>\delta ~\mbox{for all}~ n\in \mathbb{N}. \end{equation}  \end{definition}

\begin{definition} \cite{GB,G3} A sequence of real numbers $\{x_{n}\}_{n\in \mathbb{N}}$ is said to converge in a weighted manner to $\ell$ with respect to the weighted sequence $t=\{t_{n}\}_{n\in \mathbb{N}}$ if, for every $\varepsilon>0,$ there exists a positive integer $k$ such that 
	\begin{equation} \label{Eq(2)}
		t_n |x_n - \ell| < \varepsilon ~\mbox{holds for all}~ n \geq k. \end{equation}	
	This relationship is expressed as $x_n \xrightarrow{W^t} \ell.$ \end{definition}

\begin{remark} If the weighted sequence  $t = \{ t_n \}_{n \in \mathbb{N}}$ does not satisfy the required condition of Equation  \ref{Eq(1)}, i.e.,	\begin{equation}
		t_{n} > \delta ~\mbox{for all}~  n \in \mathbb{N},  ~\mbox{where}~ \delta > 0   ~\mbox{is a constant},\end{equation}
	a significant problem may occur.  For example, let $\lim t_n=0.$  From Equation \ref{Eq(2)}, we can conclude that any bounded sequence $x = \{x_n\}_{n \in \mathbb{N}}$ is weighted convergent to any real number. However, this conclusion lacks meaningful interpretation.	\end{remark}

If a sequence is convergent in the usual sense, it may not be a weighted convergent sequence. The following example gives the required answer.

\begin{example}\label{Example 1.1} Let $\{x_{n}\}_{n\in \mathbb{N}}$ be a sequence of real numbers such that $x_{n}=\frac{1}{n} $ for all $n \in \mathbb{N},$ and a weighted sequence be $t_{n}=n^2$ for all $n \in \mathbb{N}.$
	Then, $t_{n}|x_{n}-0|=n\rightarrow \infty \ \mbox{as}\ n \rightarrow \infty.$ Therefore, $x_{n}\not\xrightarrow {W^t}0.$  In conclusion, we say that the usual  and weighted convergence are two entirely different notions. \end{example}


\begin{definition}\label{Def 1.11} \cite{GM1} Let $\mathcal{I}$ be an admissible ideal of $\mathbb{N}$ and $t=\{t_{n}\}_{n\in \mathbb{N}}$ be a \emph{weighted sequence} of real numbers. A sequence $x=\{x_{n}\}_{n\in \mathbb{N}}$ in a \emph{locally solid Riesz space} $(L,\tau)$ is said to be \emph{rough weighted $\mathcal{I}_\tau$-convergent} to $x_*\in L$ w.r.t the roughness of degree $\textbf{V}_\theta$ (where $\textbf{V}_\theta$ is a $\tau$-neighbourhood of $\theta$), denoted by $x_{n}\xrightarrow[\textbf{V}_\theta]{\text{W}\mathcal{I}} x_*,$ if for every $\tau$-neighbourhood $\textbf{U}_\theta$ of $\theta,$   the following expression holds,
	$$\{n\in \mathbb{N}: t_n(x_n-x_*)\notin \textbf{V}_\theta+\textbf{U}_\theta\}\in \mathcal{I}~\text{or}, ~\{n\in \mathbb{N}: t_n(x_n-x_*)\in \textbf{V}_\theta+\textbf{U}_\theta\}\in \mathcal{F(\mathcal{I})}.$$	
	We shall write $\text{W}\mathcal{I}_{\tau}-\text{LIM}^{\textbf{V}_\theta}x=\{x_*\in L: x_{n}\xrightarrow[\textbf{V}_\theta]{\text{W}\mathcal{I}} x_*\}$ to denote the set of all \emph{rough weighted $\mathcal{I}_\tau$-limit points} of the sequence $x=\{x_{n}\}_{n\in \mathbb{N}}$ with degree of roughness $\textbf{V}_\theta.$
\end{definition}




\begin{definition} \label{Def 1.12} \cite{GM1} Let $\mathcal{I}$ be an admissible ideal of $\mathbb{N}$ and $t=\{t_{n}\}_{n\in \mathbb{N}}$ be a \emph{weighted sequence} of real numbers. An element $\zeta\in L$ is called \emph{weighted $\mathcal{I}_{\tau}$-cluster point} of a sequence $x=\{x_n\}_{n\in \mathbb{N}}$ in $L$ if for every $\tau$-neighbourhood $\textbf{U}_\theta$ of $\theta,$ the set
	$$\{n\in \mathbb{N}: t_n(x_n-\zeta)\in \textbf{U}_\theta\}\notin \mathcal{I}.$$
	The set of all \emph{weighted $\mathcal{I}_\tau$-cluster points} of the sequence $x=\{x_n\}_{n\in \mathbb{N}}$ is denoted by $\text{W}_{\mathcal{I}}\Gamma_x^\tau.$	
\end{definition}

\begin{remark} In the Definition \ref{weight sequence} \cite{GB, GM1, GM2}, a positive real number $\delta$ is defined in the context of \emph{weighted sequence} $t=\{t_{n}\}_{n\in \mathbb{N}},$ such that $t_n > \delta$ for every $n \in \mathbb{N}.$ We do not follow the usual convention and propose a new hypothesis in the Definition \ref{Def 2.1}, as follows: $$ \text{there exists a} ~\delta>0~\mbox{such that}~ |t_{\alpha}|>\delta~\mbox{for all}~\alpha \in D, ~\mbox{where}~ D ~\mbox{is a directed set}.$$ \end{remark}


\begin{definition} \label{Def 2.1} Let $(D,\geq)$ be a \emph{directed set}. A net $t=\{t_{\alpha}\}_{\alpha\in D}$ in $\mathbb{R}$ is said to be a \emph{weighted net}, if there exists a positive real number $\delta$ such that $|t_{\alpha}|>\delta$ for all $\alpha \in D.$  \end{definition}

Continuing from our previous discussion, we now present two important definitions.

\begin{definition} \label{Def 2.2}Let $\mathcal{I}$ be a non-trivial ideal of the directed set $D.$ A \emph{weighted net} $t=\{t_{\alpha}\}_{\alpha\in D}$ is said to be \emph{$\mathcal{I}$-bounded}, if there exists a positive real number $M$ such that $$\{\alpha\in D: |t_{\alpha}|>M\}\in \mathcal{I}~\mbox{or},~\{\alpha\in D: |t_{\alpha}|\leq M\}\in \mathcal{F}(\mathcal{I}).$$   \end{definition}


\begin{definition} Let $(L,\tau)$ be a \emph{locally solid Riesz space} and $\mathcal{I}$ be a non-trivial ideal of the directed set $D.$  A net $s=\{s_{\alpha}\}_{\alpha\in D}$ in $L$ is said to be \emph{weighted  $\mathcal{I}_\tau$-bounded} corresponding to the \emph{weighted net} $t=\{t_{\alpha}\}_{\alpha\in D}$ if for every $\tau$-neighbourhood ${\textbf{U}_\theta}$ of $\theta$ there exists a $\lambda>0$ such that  $$\{\alpha\in D: \lambda t_\alpha s_\alpha \notin {\textbf{U}_\theta}\}\in \mathcal{I}~\mbox{or,}~\{\alpha\in D: \lambda t_\alpha s_\alpha \in {\textbf{U}_\theta}\}\in \mathcal{F}(\mathcal{I}).$$ \end{definition}



In this paper, we focus on the concepts of ~`\emph{rough weighted ideal-$\tau$-convergence of a net}' (or, more briefly, ~`\emph{rough weighted $\mathcal{I}_\tau$-convergence of a net}') and  ~`\emph{rough weighted ideal-$\tau$-cluster points set of a net}' (or, more briefly, ~`\emph{rough weighted $\mathcal{I}_\tau$-cluster points set of a net}') within the framework of \emph{locally solid Riesz spaces}.  This article presents several significant examples and results that emphasize the importance of the concepts of  ~`\emph{rough weighted $\mathcal{I}_\tau$-convergence}' and ~`\emph{rough weighted $\mathcal{I}_\tau$-cluster points}' of a net in locally solid Riesz spaces. We present two consecutive strong examples, specifically Examples \ref{Example 2} (i) and \ref{Example 2} (ii), to demonstrate that the set of ~`\emph{rough weighted $\mathcal{I}_\tau$-cluster points of a net}' in a locally solid Riesz space may not be closed. This is why we are investigating this area. Ultimately, these new concepts allow us to establish a necessary and sufficient condition for determining when an ideal is maximal. Additionally, we will precisely characterize a closed $\tau$-neighbourhood of $\theta,$ in terms of the ~`\emph{rough weighted $\mathcal{I}_\tau$-cluster points set of a net}'.




\section{Main Results} 

Let us begin this section by asking how to define `\emph{rough weighted $\mathcal{I}_\tau$-convergence}' and 
`\emph{rough weighted $\mathcal{I}_{\tau}$-cluster point}' of a net in locally solid Riesz spaces. From the Equations \ref{Eq1.1}, \ref{Eq1.2}, we redefine the concept of `\emph{rough statistical convergence}' and `\emph{rough statistical cluster point}' of a  sequence  in a normed linear space $X$ using neighbourhoods of $\theta$, as follows:

\begin{definition}	Let $X$ be a normed linear space, $B(\theta,M)=\{x\in X : ~\parallel x - \theta \parallel <M\}$ where $M>0$ and $$B(\theta,M_1)+ B(\theta,M_2)=\{x_1+x_2 \in X : x_1 \in B(\theta,M_1) ~\mbox{and}~ x_2 \in B(\theta,M_2)\}.$$	
\begin{description}	
	\item[(i)]  A sequence  $x=\{x_n\}_{n\in \mathbb{N}}$ in $X$  is said to be rough statistically convergent to $x_*$ w.r.t the roughness degree $r> 0,$ if for any $\varepsilon > 0,$ such that 
	\begin{equation}\label{Eq(N1.1)}
		\displaystyle{\lim_{n\rightarrow \infty}}~\frac{1}{n}|\{k\leq {n}: x_k-x_* \not\in  B(\theta,r)+ B(\theta,\varepsilon)\}|=0,
	\end{equation} 
\item[(ii)] An element $\zeta$ in a norm space $X$ is called a  rough statistical cluster point w.r.t the roughness degree $r> 0$ of sequence  $x=\{x_n\}_{n\in \mathbb{N}}$ if for any $\varepsilon > 0,$ such that 
\begin{equation}\label{Eq(N1.2)}
	\displaystyle{\lim_{n\rightarrow \infty}}~\frac{1}{n}|\{k\leq {n}: x_k-\zeta \in B(\theta,r)+ B(\theta,\varepsilon)\}|\neq 0.	
\end{equation}
\end{description}
\end{definition}


According to the Equations \ref{Eq(N1.1)} and \ref{Eq(N1.2)}, we are now able to introduce the definitions of ~`\emph{rough weighted $\mathcal{I}_\tau$-convergence}' and ~`\emph{rough weighted $\mathcal{I}_{\tau}$-cluster point}' of a net in locally solid Riesz spaces as follows:


\begin{definition}\label{Def 2.4} Let $\mathcal{I}$ be a non-trivial ideal of a directed set $D$ and $t=\{t_{\alpha}\}_{\alpha\in D}$ be a \emph{weighted net}. A net $s=\{s_{\alpha}\}_{\alpha\in D}$ in a \emph{locally solid Riesz space} $(L,\tau)$ is said to be \emph{rough weighted $\mathcal{I}_\tau$-convergent} to $s_0\in L$ w.r.t the roughness of degree $\textbf{V}_\theta$ (where $\textbf{V}_\theta$ is a $\tau$-neighbourhood of $\theta$), denoted by $s_{\alpha}\xrightarrow[\textbf{V}_\theta]{\text{W}\mathcal{I}} s_0,$ if for every $\tau$-neighbourhood $\textbf{U}_\theta$ of $\theta,$   the following expression holds,
	$$\{\alpha\in D: t_\alpha(s_\alpha-s_0)\notin \textbf{V}_\theta+\textbf{U}_\theta\}\in \mathcal{I}~\text{or}, ~\{\alpha\in D: t_\alpha(s_\alpha-s_0)\in \textbf{V}_\theta+\textbf{U}_\theta\}\in \mathcal{F(\mathcal{I})}.$$	
	We shall write $\text{W}\mathcal{I}_{\tau}-\text{LIM}^{\textbf{V}_\theta}_t s=\{s_0\in L: s_{\alpha}\xrightarrow[\textbf{V}_\theta]{\text{W}\mathcal{I}} s_0,\}$ to denote the set of all \emph{rough weighted $\mathcal{I}_\tau$-limit points} of the net $s=\{s_{\alpha}\}_{\alpha\in D}$ with degree of roughness $\textbf{V}_\theta.$
\end{definition}



\begin{definition}\label{Def 2.5} Let $\mathcal{I}$ be an non-trivial ideal of a directed set $D,$ $L$ be a \emph{locally solid Riesz space} and $t=\{t_{\alpha}\}_{\alpha\in D}$ be a \emph{weighted net}. A point $\zeta\in L$ is called a \emph{rough weighted $\mathcal{I}_\tau$-cluster point of the net} $s=\{s_{\alpha}\}_{\alpha\in D}$ w.r.t the roughness of degree $\textbf{V}_\theta$ (where $\textbf{V}_\theta$ is a $\tau$-neighbourhood of $\theta$), if for every $\tau$-neighbourhood $\textbf{U}_\theta$ of $\theta,$
		$$ \{\alpha \in D: t_\alpha(s_\alpha-\zeta)\in \textbf{V}_\theta+\textbf{U}_\theta\}\notin \mathcal{I}.$$	
	The set of all \emph{rough weighted $\mathcal{I}_\tau$-cluster points of a net} $s=\{s_{\alpha}\}_{\alpha\in D}$ is denoted by $\mathcal{R}_{t}^{\textbf{V}_\theta}-{\mathcal{I}}(\Gamma_{s}^{\tau}).$ 
\end{definition}

\begin{remark} In a \emph{locally solid Riesz space} $(L,\tau)$ the following condition holds:
	
	 \begin{description}
		\item[(i)] 	Let $\mathcal{I}$ be a non-trivial ideal of a directed set $D,$  $t=\{t_{\alpha}\}_{\alpha\in D}$ be a \emph{weighted net} and  $s=\{s_{\alpha}\}_{\alpha\in D}$ be any net in $L.$ Then $\text{W}\mathcal{I}_{\tau}-\text{LIM}^{\textbf{V}_\theta}_t s \subseteq \mathcal{R}_{t}^{\textbf{V}_\theta}-{\mathcal{I}}(\Gamma_{s}^{\tau}).$
		
		\item[(ii)] Let $\mathcal{I}$ be a non-trivial ideal of $\mathbb{N},$ $t=\{t_{n}\}_{n\in \mathbb{N}}$ be a \emph{weighted sequence} and $x=\{x_{n}\}_{n\in \mathbb{N}}$ be a sequence in $L.$ 
		From the Definition \ref{Def 1.11}, \ref{Def 1.12} and \ref{Def 2.5} we get		
		  $$\text{W}_{\mathcal{I}}\Gamma_x^\tau\subseteq \mathcal{R}_{t}^{\textbf{V}_\theta}-{\mathcal{I}}(\Gamma_{x}^{\tau})~\mbox{and}~ \text{W}\mathcal{I}_{\tau}-\text{LIM}^{\textbf{V}_\theta}x \subseteq \mathcal{R}_{t}^{\textbf{V}_\theta}-{\mathcal{I}}(\Gamma_{x}^{\tau}).$$     
	\end{description}	
		\end{remark}



The following example gives the disparity between the notions of ~`\emph{rough weighted $\mathcal{I}_{\tau}$-convergence nets}' and ~`\emph{rough weighted $\mathcal{I}_{\tau}$-cluster point of a nets}' on \emph{locally solid Riesz spaces} i.e.,
 $$\left(\mathcal{R}_{t}^{\textbf{V}_\theta}-{\mathcal{I}}(\Gamma_{s}^{\tau}) \right) \setminus \left( \text{W}\mathcal{I}_{\tau}-\text{LIM}^{\textbf{V}_\theta}_t s \right) \neq \varnothing.$$


\begin{example} \label{Example 1} (i) Let us consider the \emph{locally solid Riesz space} $(\mathbb{R}^2,  ||\cdot||)$ with the Euclidean norm $||\cdot||$ and coordinate-wise ordering. The family $\mathcal{N}_{sol}$ of all $\textbf{U}(r)$ is defined as $\textbf{U}_\theta(r)=\{x\in \mathbb{R}^2 :||x||<r\}=\{(x,y)\in \mathbb{R}^2: x^2+y^2<r^2\}$  where $r>0,$ is a base at $\theta=(0,0).$  It is known that $\mathcal{N}_{sol}$ is itself a directed set $D$ with the binary relation $\geq$ defined by $\textbf{U}_\theta\geq \textbf{V}_\theta$ if $\textbf{U}_\theta\supseteq \textbf{V}_\theta$ for any $\textbf{U}_\theta,\textbf{V}_\theta\in D.$	Let $\textbf{W}_\theta \in D=\mathcal{N}_{sol},$ $D_{\textbf{W}_\theta}=\{\textbf{Z}_\theta\in D: \textbf{Z}_\theta \geq \textbf{W}_\theta\},$ then the collection 
	$$\mathcal{F}_0=\{A\subset D: ~\mbox{there exists a set}~D_{\textbf{W}_\theta} ~\mbox{with}~ A\supset D_{\textbf{W}_\theta}\}$$ forms a filter in $D.$ Let $\mathcal{I}_0=\{A\subset D: A^c\in \mathcal{F}_0\}.$ Then $\mathcal{I}_0$ is a non-trivial ideal of $D.$ Also $G=\{\textbf{U}_\theta(r): r\in \mathbb{R}^+\setminus \mathbb{Q}^+\}$ neither belongs to $\mathcal{I}_0$ nor to $\mathcal{F}_0.$ 
	
	 Next we choose a net $s=\{s_{\textbf{U}_\theta}\}_{\textbf{U}_\theta\in D},$ a weighted net $t=\{t_{\textbf{U}_\theta}\}_{\textbf{U}_\theta\in D}$ and a $\tau$-neighbourhood of $\theta,$  say $\textbf{V}_\theta,$ in the following manner:	
	 \begin{equation*}
	s_{\textbf{U}_\theta(r)}\in \begin{cases} (1,0)+\{\overline{\textbf{U}_\theta(1)}\setminus \textbf{U}_\theta(1)\}  ~\mbox {if} ~\textbf{U}_\theta(r)\in G, \\
		(4,0)+\{\overline{\textbf{U}_\theta(4)}\setminus \textbf{U}_\theta(4)\}  ~ \mbox {if} ~  \textbf{U}_\theta(r)\in D\setminus G,
	\end{cases}
		~\mbox{and}~ ~
		t_{\textbf{U}_\theta(r)}=\begin{cases} \frac{e^r}{1+e^r} & \mbox {if}~  r\in G,\\
			r+1 & \mbox {if}~  r\in D\setminus G,\\
		\end{cases}
	\end{equation*}
	
	 \begin{equation*}
	\mbox{i.e.,}~	s_{\textbf{U}_\theta(r)}\in \begin{cases} (1,0)+\{(x,y)\in \mathbb{R}^2: x^2+y^2=1\}  ~\mbox {if} ~\textbf{U}_\theta(r)\in G, \\
			(4,0)+\{(x,y)\in \mathbb{R}^2: x^2+y^2=4^2\}  ~ \mbox {if} ~  \textbf{U}_\theta(r)\in D\setminus G,
		\end{cases}	
	\end{equation*}	
	 and $\textbf{V}_\theta=\overline{\textbf{U}_\theta(1)}=\{(x,y)\in \mathbb{R}^2: x^2+y^2\leq 1\} .$  Then	$$\{\textbf{U}_\theta (r)\in D: t_{\textbf{U}_\theta(r)}(s_{\textbf{U}_\theta(r)}-(1,0))\in \textbf{V}_\theta+\textbf{W}_\theta\}=G$$ neither belongs to $\mathcal{I}_0$ nor to $\mathcal{F}_0,$ for any $\tau$-neighbourhood $\textbf{W}_\theta$ of $\theta$ where $\textbf{W}_\theta \subset \textbf{V}_\theta.$ This implies   $(1,0) \in \mathcal{R}_{t}^{\textbf{V}_\theta}-{\mathcal{I}_0}(\Gamma_{s}^{\tau}),$ whereas the net $s=\{s_{\textbf{U}_\theta}\}_{\textbf{U}_\theta\in D}$ cannot \emph{rough weighted $\mathcal{I}_{\tau}$-converge} to $(1,0)$  w.r.t the roughness of degree $\textbf{V}_\theta.$	 
	 $$(1,0) \in \left(\mathcal{R}_{t}^{\textbf{V}_\theta}-{\mathcal{I}}(\Gamma_{s}^{\tau}) \right) \setminus \left( \text{W}\mathcal{I}_{\tau}-\text{LIM}^{\textbf{V}_\theta}_t s \right) \neq \varnothing.$$

 	(ii) Let $\{A_n : n\in \mathbb{N}\}$ be a partition of the directed set $D=\mathbb{R}^+$ into infinite sets as in Example \ref{Example 2} (i). Also, $\textbf{V}_\theta=(-\frac{1}{2},\frac{1}{2})$
 and  ${\textbf{U}_\theta} \subset \textbf{V}_\theta.$ 	 For each $n\in \mathbb{N}$ and $\alpha \in A_n,$ define $s_\alpha = n$ and $t_\alpha =n.$ Let $\mathcal{I}\subset \mathbb{P}(\mathbb{R}^+)$ be the ideal of sets $I$ satisfying that $I \cap A_n$ is finite for all $n\in \mathbb{N}.$  Then the set	$\{\alpha \in \mathbb{R}^+ : t_\alpha (s_\alpha- 1) = 0 \in {\textbf{V}_\theta}+ {\textbf{U}_\theta}  \}=A_1$ neither belongs to $\mathcal{I}$ nor to $\mathcal{F}(\mathcal{I}).$ This implies that  
 $$1\in \left(\mathcal{R}_{t}^{\textbf{V}_\theta}-{\mathcal{I}}(\Gamma_{s}^{\tau}) \right) \setminus \left( \text{W}\mathcal{I}_{\tau}-\text{LIM}^{\textbf{V}_\theta}_t s \right) \neq \varnothing.$$
   \end{example}

We now prove the expected relationship between the sets $ \text{W}{\mathcal{I}_{\tau}}-\text{LIM}^{\textbf{V}_\theta}_t s$ and $\mathcal{R}_{t}^{\textbf{V}_\theta}-{\mathcal{I}}(\Gamma_{s}^{\tau}).$


\begin{theorem}   Let $\mathcal{I}^D$ be an $D$-admissible ideal of a directed set $D,$ $t=\{t_{\alpha}\}_{\alpha\in D}$ be a \emph{weighted net}, $s=\{s_{\alpha}\}_{\alpha\in D}$ be a net in a \emph{locally solid Riesz space} $L$ and $\textbf{V}_\theta$ be a $\tau$-neighbourhood of $\theta.$ If $\zeta \in \mathcal{R}_{t}^{\textbf{V}_\theta}-{\mathcal{I}^D}(\Gamma_{s}^{\tau}),$ then there exists a subnet $\omega=\{\omega_{e}\}_{e\in E}$ of the net $\{t_\alpha (s_\alpha-\zeta)\}_{\alpha\in D}$ such that
	 $\theta \in  \text{W}{\mathcal{I}^{E}_{\tau}}-\text{LIM}^{\textbf{W}_\theta}_1 \omega,$ where $\mathcal{I}^{E}$ is a $E$-admissible ideal of $E$ and $\textbf{W}_\theta$ is $\tau$-neighbourhood of $\theta$ such that $\textbf{V}_\theta \subseteq \textbf{W}_\theta.$
	\end{theorem}

\begin{proof} 	 Since  $\zeta \in \mathcal{R}_{t}^{\textbf{V}_\theta}-{\mathcal{I}^D}(\Gamma_{s}^{\tau}),$ then for every $\tau$-neighbourhood $\textbf{U}_\theta$ of $\theta,$ the set
	$$A=\{\alpha \in D: t_\alpha (s_\alpha-\zeta)\in \textbf{V}_\theta+\textbf{U}_\theta\}\not\in \mathcal{I}^D$$ 
	$$\implies ~D \setminus A=\{\alpha \in D: t_\alpha (s_\alpha-\zeta) \notin \textbf{V}_\theta+\textbf{U}_\theta\}\not\in \mathcal{F}(\mathcal{I}^D).$$	
	Let $\alpha_0\in D.$ Then $D_{\alpha_0}=\{\alpha \in D: \alpha \geq \alpha_0\}\in \mathcal{F}(\mathcal{I}^D) $ and $D_{\alpha_0}\not\subseteq D\setminus A$ (otherwise $D\setminus A \in \mathcal{F}(\mathcal{I}^D)$ a contradiction). There is an element $\alpha_* \in D_{\alpha_0}$ such that $\alpha_* \not\in D\setminus A.$ 	
	This implies $\alpha_* \in  A \subset D,$  $t_{\alpha_*} (s_{\alpha_*}-\zeta)\in \textbf{V}_\theta+\textbf{U}_\theta$ and $\alpha_* \geq \alpha_0.$ So $\alpha_*$ depending on $\alpha_0$ and $\textbf{U}_\theta.$ We construct a directed set $E,$ a subnet $\omega=\{\omega_{e}\}_{e\in E}$ of the net $\{t_\alpha (s_\alpha-\zeta)\}_{\alpha\in D}$ as follows:
	$$E=\{(\alpha, \textbf{U}_\theta)\in D\times \mathcal{N}_{\theta}: t_{\alpha} (s_{\alpha}-\zeta)\in \textbf{V}_\theta+\textbf{U}_\theta \},$$ 
	where $\mathcal{N}_{\theta}$ denote the family of all $\tau$-neighbourhood of $\theta.$  For $(\alpha, \textbf{U}_\theta),$  $(\beta, \textbf{U}^*_\theta) \in E$ define a relation $(\beta, \textbf{U}^*_\theta) \geq (\alpha, \textbf{U}_\theta)$ mean $\beta\geq \alpha$ and $\textbf{U}^*_\theta \subseteq \textbf{U}_\theta.$ Then $E$ is a directed set.	Now we choose two mappings $i: E\rightarrow D$ and $\omega : E \rightarrow L$ as follows:
	$$i((\alpha_0,\textbf{U}_\theta))=\alpha_* ~\mbox{and}~ \omega((\alpha_0,\textbf{U}_\theta))=t_{\alpha_*} (s_{\alpha_*}-\zeta)~\mbox{for all}~(\alpha_0,\textbf{U}_\theta)\in E.$$
	
	Let $\textbf{U}^*_\theta$ be a $\tau$-neighbourhood of $\theta$ and $\beta \in D.$ Then there exists an element $\beta_{*} \in D$ depends on $\textbf{U}^*_\theta$ and $\beta$ with $\beta_* \geq \beta$ such that 
	$$\omega((\beta,\textbf{U}^*_\theta))=t_{\beta_{*}} (s_{\beta_{*}}-\zeta) \in \textbf{V}_\theta+\textbf{U}^*_\theta \subseteq \textbf{V}_\theta+\textbf{U}_\theta \subseteq \textbf{W}_\theta+\textbf{U}_\theta, $$
	for all $(\beta,\textbf{U}^*_\theta) \geq (\alpha_0,\textbf{U}_\theta),$ $\beta\geq \alpha_0$ and $\textbf{U}^*_\theta \subseteq \textbf{U}_\theta.$ Therefore	
	$$E_{(\alpha_0,\textbf{U}_\theta)}\subseteq \{(\beta,\textbf{U}^*_\theta) \in E: \omega((\beta,\textbf{U}^*_\theta)) \in \textbf{W}_\theta+\textbf{U}_\theta\}\in \mathcal{F}(\mathcal{I}^E).$$
 	Therefore this is the result.	
\end{proof}


\begin{theorem} \label{subnet2}  Let $D$ be a directed set, $t=\{t_{\alpha}\}_{\alpha\in D}$ be a \emph{weighted net}, $s=\{s_{\alpha}\}_{\alpha\in D}$ be a net in a \emph{locally solid Riesz space} $L,$  $\zeta \in L,$ $\textbf{W}_\theta$ be a $\tau$-neighbourhood of $\theta$  and $\omega=\{\omega_{e}\}_{e\in E}$ be a subnet of the net $\{t_\alpha (s_\alpha-\zeta)\}_{\alpha\in D}.$  If 	$\theta \in  \text{W}{\mathcal{I}^{E}_{\tau}}-\text{LIM}^{\textbf{W}_\theta}_1 \omega,$ where $\mathcal{I}^{E}$ is a $E$-admissible ideal of $E$ then $\zeta \in \mathcal{R}_{t}^{\textbf{V}_\theta}-{\mathcal{I}^D}(\Gamma_{s}^{\tau}),$ for some $D$-admissible ideal $\mathcal{I}^D$
 and $\textbf{W}_\theta \subseteq \textbf{V}_\theta.$
\end{theorem}

\begin{proof} Since $\omega=\{\omega_{e}\}_{e\in E}$ is a subnet of the net $\{t_\alpha (s_\alpha-\zeta)\}_{\alpha\in D}$ there exists a mapping $i: E\rightarrow D$ such that $\omega(e)=Soi(e)$ for all $e\in E,$ where $S_\alpha=t_\alpha(s_\alpha-\zeta)$ for all $\alpha\in D.$ As in Remark \ref{new ideal}, the collection $\mathcal{I}^D=\{A\subset D: i^{-1}(A)\in \mathcal{I}_E\}$ is a $D$-admissible ideal of $D.$ We assume by the way of contradiction  $\zeta \not\in \mathcal{R}_{t}^{\textbf{V}_\theta}-{\mathcal{I}^D}(\Gamma_{s}^{\tau}),$ then the set 
		$$A=\{\alpha \in D: t_\alpha (s_\alpha-\zeta)\in \textbf{V}_\theta+\textbf{U}_\theta\}\in \mathcal{I}^D~\mbox{and}~i^{-1}(A)\in \mathcal{I}^E.$$

Since $\theta \in  \text{W}{\mathcal{I}^{E}_{\tau}}-\text{LIM}^{\textbf{W}_\theta}_1 \omega,$ then the set
$$C=\{e \in E:  \omega_e\in \textbf{W}_\theta+\textbf{U}_\theta\}\in \mathcal{F}(\mathcal{I}^E).$$

Let $e\in C.$ Then $$\omega_e\in \textbf{W}_\theta+\textbf{U}_\theta~\implies~Soi(e) \in \textbf{W}_\theta+\textbf{U}_\theta ~\implies~S(i(e)) \in \textbf{W}_\theta+\textbf{U}_\theta\subseteq \textbf{V}_\theta+\textbf{U}_\theta,$$
$\implies~ i(e)\in A ~\implies~e\in i^{-1}(A).$
Therefore $C\subseteq i^{-1}(A).$ Since $C\in \mathcal{F}(\mathcal{I}^E)$ so $i^{-1}(A) \in \mathcal{F}(\mathcal{I}^E),$ which is a contradiction since $i^{-1}(A)\in \mathcal{I}^E.$ \end{proof}

As an immediate consequence of the above two theorems we get following remark.

\begin{remark}  Let $D$ be a directed set, $t=\{t_{\alpha}\}_{\alpha\in D}$ be a \emph{weighted net}, $s=\{s_{\alpha}\}_{\alpha\in D}$ be a net in a \emph{locally solid Riesz space} $L,$ $\textbf{V}_\theta$ be a $\tau$-neighbourhood of $\theta$ and $\zeta \in L.$  Then $\zeta \in \mathcal{R}_{t}^{\textbf{V}_\theta}-{\mathcal{I}^D}(\Gamma_{s}^{\tau})$ for some $D$-admissible ideal $\mathcal{I}^D$ of a directed set $D$ {\bf iff} there exists a subnet $\omega=\{\omega_{e}\}_{e\in E}$ of the net $\{t_\alpha (s_\alpha-\zeta)\}_{\alpha\in D}$ such that
	$\theta \in  \text{W}{\mathcal{I}^{E}_{\tau}}-\text{LIM}^{\textbf{V}_\theta}_1 \omega,$ where $\mathcal{I}^{E}$ is a $E$-admissible ideal of $E.$ 
\end{remark}

Recall an important result of Lahiri et al. \cite[Theorem 9]{Lahiri}  ~`{\it In a topological space the $\mathcal{I}$-cluster points set of a net is closed }'. While reinventing this theorem in terms of ~`\emph{rough weighted $\mathcal{I}_\tau$-cluster points set of a net} in a \emph{locally solid Riesz space}', we assert that \emph{rough weighted $\mathcal{I}_\tau$-cluster points set of a net} in a \emph{locally solid Riesz space} is not necessarily closed. Our affirmation is emphatically demonstrated in the following example.


\begin{example}\label{Example 2} (i) Consider the \emph{locally solid Riesz space} $(\mathbb{R}^2,||\cdot||)$ in conjunction with the Euclidean norm $||\cdot||$ and coordinate-wise ordering. The directed set $D$ of all ${\textbf{U}_\theta}(r)$ is defined as ${\textbf{U}_\theta}(r)=\{x\in \mathbb{R}^2 :||x||<r\}=\{(x,y)\in \mathbb{R}^2: x^2+y^2<r^2\}$ where $r>0$ with the binary relation $\geq$ defined by ${\textbf{U}_\theta}\geq {\textbf{V}_\theta}$ if ${\textbf{U}_\theta}\subseteq {\textbf{V}_\theta}$ for any ${\textbf{U}_\theta},{\textbf{V}_\theta}\in D.$ Consider a new decomposition of $\mathbb{R}^+,$  
		
		$$B_0=(0,1)\cup \{\displaystyle{\bigcup_{n\in O}}~[n,n+1]\}~\mbox{where} ~O=\{1,3,5,\ldots\},$$
		
		$$B_\ell=\displaystyle{\bigcup_{n\in O}}~( 2^{\ell} n, 2^{\ell}n+1)~\mbox{where} ~ \ell\in \mathbb{N}=\{1,2,3,\ldots\},$$

		$$\implies~\mathbb{R}^+=\displaystyle{\bigcup_{\ell \in \mathbb{N}_0}}~B_\ell~\mbox{where}~~B_i \cap B_j=\varnothing ~\mbox{for all}~ i,j\in \mathbb{N}_0=\{0,1,2,3,\ldots\} ~\mbox{and}~i\neq j.$$

		Let  $\mathcal{I}=\{A\subset D: A ~\mbox{is  countable}\}$ be a proper non-trivial ideal of $D$ and $\textbf{V}_\theta=\overline{\textbf{U}_\theta(1)}=\{(x,y)\in \mathbb{R}^2: x^2+y^2\leq 1\} .$	Now we select a net $s=\{s_{\alpha}\}_{\alpha\in D}$ and a weighted net $t=\{t_{\alpha}\}_{\alpha\in D}$ in the following manner:
		\begin{equation*}
			s_{{\textbf{U}_\theta}(r)}=\begin{cases} (2+r,1)  &\mbox {if} ~ r\in B_0, \\
				(\frac{1}{2^\ell}+\frac{1}{r}+\frac{1}{r^2},1) & \mbox {if} ~ r\in B_{\ell}, ~\ell\in \mathbb{N},
			\end{cases}
			~\mbox{and}~ ~~		
			t_{{\textbf{U}_\theta}(r)}= \begin{cases} 2  & \mbox {if} ~  r\in B_0, \\		
				r+1 & \mbox{if}~  r\in B_{\ell}, ~\ell\in \mathbb{N}.
			\end{cases}
		\end{equation*}
		
		Let ${\textbf{U}_\theta}$ be a $\tau$-neighbourhood of $\theta,$ there exists some ${\textbf{U}_\theta}(\varepsilon)\in \mathcal{N}_{sol},$ $0<\varepsilon<\frac{1}{2}$ such that ${\textbf{U}_\theta}(\varepsilon)=\{(x,y)\in \mathbb{R}^2: x^2+y^2<\varepsilon^2\}\subset {\textbf{U}_\theta}.$ Moreover for each $m \in \mathbb{N},$ we derive that	
		\begin{equation*}	
			{t_{{\textbf{U}_\theta}(r)}}\left(s_{{\textbf{U}_\theta}(r)}-(\frac{1}{2^m},1)\right) =\begin{cases} \left(4+2r-\frac{1}{2^{m-1}},0\right) & \mbox{if}~ r\in B_0, \\	
				\left( (r+1)(\frac{1}{2^{\ell}}-\frac{1}{2^m}) +1+\frac{2}{r}+\frac{1}{r^2},0 \right)  &\mbox{if}~ r\in B_{\ell}, ~\ell\in \mathbb{N}.
			\end{cases} 
		\end{equation*}

		
		Let $r\in B_{\ell}, ~\ell\in \mathbb{N}$ and $m=\ell$ then $$	{t_{{\textbf{U}_\theta}(r)}}\left(s_{{\textbf{U}_\theta}(r)}-(\frac{1}{2^m},1)\right) =
		\left(1+\frac{2}{r}+\frac{1}{r^2},0 \right)\in \textbf{V}_\theta+ {\textbf{U}_\theta}(\varepsilon),$$	 
		$$\implies~\{{\textbf{U}_\theta}(r)\in D: {t_{{\textbf{U}_\theta}(r)}}\left(s_{{\textbf{U}_\theta}(r)}-(\frac{1}{2^m},1)\right)\in \textbf{V}_\theta+{\textbf{U}_\theta}(\varepsilon)\}\notin \mathcal{I} .$$	
		This conveys $(\frac{1}{2^m},1)\in \mathcal{R}_{t}^{\textbf{V}_\theta}-{\mathcal{I}}(\Gamma_{s}^{\tau})$ for all $m\in \mathbb{N}.$ For $r\in B_{\ell}$ where $\ell \in \mathbb{N}$ there exists an odd positive integer $n_0 \in O$ such that $r \in (2^{\ell}n_0, 2^{\ell}n_0+1)$ then $r> 2^{\ell}n_0$ implies  $\frac{r+1}{2^\ell}\geq n_0\geq 1.$  Hence we get 	
		$$\{{\textbf{U}_\theta}(r)\in D: {t_{{\textbf{U}_\theta}(r)}}\left(s_{{\textbf{U}_\theta}(r)}-(0,1)\right)\in \textbf{V}_\theta+{\textbf{U}_\theta}(\varepsilon)\}$$
		\begin{equation*} 
			= \begin{cases} \{{\textbf{U}_\theta}(r)\in D, r\in B_0: \left(4+2r,0\right)\in \textbf{V}_\theta+ {\textbf{U}_\theta}(\varepsilon)\}=\{\mbox{Finite set}\} \in \mathcal{I},\\
				\{{\textbf{U}_\theta}(r)\in D,  r\in B_\ell ~\mbox{for}~ \ell \in \mathbb{N}: \left(1+\frac{r+1}{2^{\ell}}+\frac{2}{r}+\frac{1}{r^2},0 \right)\in \textbf{V}_\theta+{\textbf{U}_\theta}(\varepsilon)\}=\{\mbox{Finite set}\} \in \mathcal{I}.
			\end{cases}  	
		\end{equation*}		
		Above illustration emphasizes $(0,1)\notin \mathcal{R}_{t}^{\textbf{V}_\theta}-{\mathcal{I}}(\Gamma_{s}^{\tau})$  and so our affirmation $\mathcal{R}_{t}^{\textbf{V}_\theta}-{\mathcal{I}}(\Gamma_{s}^{\tau})$ is not a closed set is established. \\
		
	(ii) Let $\{x_n\}_{n\in \mathbb{N}} \subset (0,1)$ be a sequence converging to $0$ and $\{A_n : n\in \mathbb{N}\}$ be a partition of $D=\mathbb{R}^+$ into infinite sets as above. Also, $\textbf{V}_\theta=(-1,1)$
	and  ${\textbf{U}_\theta} \subset \textbf{V}_\theta.$ 	 For each $n\in \mathbb{N}$ and $\alpha \in A_n,$ define $s_\alpha = x_n$ and $t_\alpha =\frac{3}{x_n}.$ Let $\mathcal{I}\subset \mathbb{P}(\mathbb{R}^+)$ be the ideal of sets $I$ satisfying that $I \cap A_n$ is finite for all $n\in \mathbb{N}.$  Then, for each $n,$ the set	$\{\alpha \in \mathbb{R}^+ : t_\alpha (s_\alpha- x_n) = 0 \in {\textbf{V}_\theta}+ {\textbf{U}_\theta}  \}$ contains the set $A_n$ and so is not in $\mathcal{I}.$	
	This implies that all the point of the sequence $\{x_n\}_{n\in \mathbb{N}}$ is a rough weighted $\mathcal{I}_{\tau}$-cluster point of the net $s=\{s_{\alpha}\}_{\alpha\in D} $ i.e., $x_n\in \mathcal{R}_{t}^{\textbf{V}_\theta}-{\mathcal{I}}(\Gamma_{s}^{\tau})$ for all $n\in \mathbb{N}.$ On the other hand,  the set $\{\alpha \in \mathbb{R}^+ : t_\alpha(s_\alpha - 0) = 3 \in {\textbf{V}_\theta}+ {\textbf{U}_\theta} \}$ is empty, and so $0$ is $\bf{not}$  rough weighted $\mathcal{I}_{\tau}$-cluster point of the net $s=\{s_{\alpha}\}_{\alpha\in D}$ i.e., $0\notin \mathcal{R}_{t}^{\textbf{V}_\theta}-{\mathcal{I}}(\Gamma_{s}^{\tau}).$ Therefore, set $\mathcal{R}_{t}^{\textbf{V}_\theta}-{\mathcal{I}}(\Gamma_{s}^{\tau})$ is not closed.	
\end{example}
	

An important theorem is now presented, derived from the above example. We show the following results holds without using the \emph{locally solid Riesz spaces} are first countable, Hausdorff or,  $\textbf{V}_\theta$~ is ~$\tau$-bounded.


\begin{prop} Let $\mathcal{I}$ be a non-trivial ideal of a directed set $D$ and $(L,\tau)$ be a locally solid Riesz space. If the weighted net  $t=\{t_{\alpha}\}_{\alpha\in D}$ is $\mathcal{I}$-bounded, then the set  $\mathcal{R}_{t}^{\textbf{V}_\theta}-\mathcal{I}(\Gamma_{s}^{\tau})$ is closed. \end{prop}

\begin{proof} From the $\mathcal{I}$-boundedness of the \emph{weighted net}, there exists a positive real number $M$ such that the set $A=\{\alpha\in D: |t_\alpha|< M\}\in \mathcal{F(\mathcal{I})}.$
	If $\textbf{U}_\theta$ be any neighbourhood of $\theta,$ then there exists a $\textbf{W}_\theta\in \mathcal{N}_{sol}$ such that $\textbf{W}_\theta+ \textbf{W}_\theta \subseteq \textbf{U}_\theta.$  Let $\omega \in \overline{\mathcal{R}_{t}^{\textbf{V}_\theta}-{\mathcal{I}}(\Gamma_{s}^{\tau})}.$ Then there exists a $$\zeta \in \{\mathcal{R}_{t}^{\textbf{V}_\theta}-{\mathcal{I}}(\Gamma_{s}^{\tau})\}\cap (\omega+\frac{1}{M}\textbf{W}_\theta).$$ This implies $(\zeta-\omega)M\in \textbf{W}_\theta$ and $B=\{\alpha\in D: t_\alpha(s_\alpha-\zeta)\in \textbf{V}_\theta+\textbf{W}_\theta\}\notin \mathcal{\mathcal{I}}.$ 
	Therefore $A\cap B\notin \mathcal{I}$ (otherwise, $B\subset (\mathbb{N}\setminus A)\cup (A\cap B)\in \mathcal{I},$ a contradiction). Let $\alpha \in A\cap B$ then 
	$$t_\alpha(s_{\alpha}-\omega)=t_\alpha(s_\alpha-\zeta)+t_\alpha(\zeta-\omega)\in \textbf{V}_\theta+\textbf{W}_\theta+\textbf{W}_\theta\subseteq \textbf{V}_\theta+\textbf{U}_\theta,$$
	since $|t_\alpha(\zeta-\omega)|\leq |M(\zeta-\omega)|,$ $(\zeta-\omega)M\in \textbf{W}_\theta$ and $\textbf{W}_\theta \in \mathcal{N}_{sol}$ then $t_\alpha(\zeta-\omega)\in \textbf{W}_\theta.$ Therefore	 
	$$A\cap B \subseteq \{\alpha \in D: t_\alpha(s_\alpha-\omega)\in \textbf{V}_\theta+\textbf{U}_\theta\} \notin \mathcal{I}.$$ Hence the set  $\mathcal{R}_{t}^{\textbf{V}_\theta}-\mathcal{I}(\Gamma_{s}^{\tau})$ is closed.\end{proof}

The  following results that shows non-emptiness of the set $\mathcal{R}_{t}^{\textbf{V}_\theta}-{\mathcal{I}}(\Gamma_{s}^{\tau}).$


\begin{prop} If the nets $s=\{s_{\alpha}\}_{\alpha\in D}$ is {weighted $\mathcal{I}_\tau$-bounded} corresponding to the weighted net $t=\{t_{\alpha}\}_{\alpha\in D},$ then for every $\tau$-neighbourhood ${\textbf{U}_\theta}$ of $\theta,$ there exists a positive real number $\vartheta$  such that $\mathcal{R}_{t}^{(\vartheta\textbf{U}_\theta)}-{\mathcal{I}}(\Gamma_{s}^{\tau})\neq \varnothing.$ \end{prop}

\begin{proof} Since  $s=\{s_{\alpha}\}_{\alpha\in D}$  is \emph{weighted $\mathcal{I}_\tau$-bounded} then for every $\tau$-neighborhood $\textbf{U}_\theta$ of $\theta$ there exists some $\frac{1}{\vartheta}>0$ such that  $$\{\alpha \in D: t_\alpha (s_\alpha-\theta) \in \vartheta \textbf{U}_\theta\}=\{\alpha \in D: t_\alpha s_\alpha \in \vartheta \textbf{U}_\theta\} = \{\alpha \in D: \frac{1}{\vartheta} t_\alpha s_\alpha \in \textbf{U}_\theta\}\in \mathcal{F}(\mathcal{I}).$$
	
	 Therefore $\theta \in \mathcal{R}_{t}^{(\vartheta\textbf{U}_\theta)}-{\mathcal{I}}(\Gamma_{s}^{\tau}).$ As a result, this is what we have concluded  \end{proof}

 Next, we choose an example showing that the converse is not valid.

  
\begin{example} Consider the \emph{locally solid Riesz space}, directed set $D,$ $G$ and ideal $\mathcal{I}$ as defined in Example \ref{Example 1} (i). Also we select the net $s=\{s_{\alpha}\}_{\alpha\in D}$ and the weighted net $t=\{t_{\alpha}\}_{\alpha\in D}$ in respectivly:
	\begin{equation*}		s_{{\textbf{U}_\theta}(r)}=\begin{cases} (2,3)   & \mbox {if} ~ {\textbf{U}_\theta}(r)\in D\setminus {G}, \\
			(r,r)  & \mbox {if} ~             ~~~~~ {\textbf{U}_\theta}(r)\in {G},
		\end{cases}
		~\mbox{and}~
		t_{\textbf{U}_\theta(r)}= -(r+1)~ \mbox{for all}~ {\textbf{U}_\theta}(r)\in D.
	\end{equation*}
	If ${\textbf{U}_\theta}$ is any arbitrary $\tau$-neighbourhood of $\theta$ and $\vartheta$ is any positive real number, then $(2,3)\in \mathcal{R}_{t}^{(\vartheta\textbf{U}_\theta)}-{\mathcal{I}}(\Gamma_{s}^{\tau})$  and $$\{\textbf{U}_\theta (r) \in D: \vartheta t_{\textbf{U}_\theta(r)} s_{\textbf{U}_\theta(r)} \notin {\textbf{U}_\theta}(1)\}\notin \mathcal{I}.$$ 
	Hence the nets $s=\{s_{\alpha}\}_{\alpha \in D}$ is not \emph{weighted $\mathcal{I}_\tau$-bounded} whereas $\mathcal{R}_{t}^{(\vartheta\textbf{U}_\theta)}-{\mathcal{I}}(\Gamma_{s}^{\tau})\neq \varnothing.$ \end{example}

Now, we figure out some other results related to the non-emptiness of the set $\mathcal{R}_{t}^{\textbf{V}_\theta}-{\mathcal{I}}(\Gamma_{s}^{\tau}).$


\begin{prop}  Let $\mathcal{I}$ be a non-trivial ideal of a directed set $D,$ $t=\{t_{\alpha}\}_{\alpha\in D}$ be a \emph{weighted net} and $s=\{s_{\alpha}\}_{\alpha\in D}$ be a net in a \emph{locally solid Riesz space} $L.$  Suppose   $H_{\ell}=\{t_\alpha(s_\alpha-\ell) \in L: \alpha\in H\}$ where $\ell \in L$ and $H\in \mathcal{F(\mathcal{I})}.$ Then	
	\begin{description}
		\item[(i)] $\zeta\in \mathcal{R}_{t}^{\textbf{V}_\theta}-{\mathcal{I}}(\Gamma_{s}^{\tau})$ if $\theta \in \overline{H_{\zeta}} $ for all $H\in \mathcal{F(\mathcal{I})},$		
		\item[(ii)] $H_{\zeta} \cap \overline{\textbf{V}_\theta} \neq \varnothing$ for all $H\in \mathcal{F(\mathcal{I})},$  if $\zeta\in \mathcal{R}_{t}^{\textbf{V}_\theta}-{\mathcal{I}}(\Gamma_{s}^{\tau})$ and $L$ is Hausdorff.
	\end{description}	
\end{prop}

\begin{proof} {\bf (i)} If possible, let $\zeta$ be not a \emph{rough weighted $\mathcal{I}_\tau$-cluster point} of $s=\{s_{\alpha}\}_{\alpha\in D}.$ Then  there exists a  $\tau$-neighbourhood ${\textbf{U}_\theta}$ of $\theta$ such that $$H^{0}=\{\alpha\in D: t_\alpha(s_\alpha-\zeta)\in L\setminus (\textbf{V}_\theta+\textbf{U}_\theta)\}\in \mathcal{F(\mathcal{I})},$$  $$H_{\zeta}^0=\{t_\alpha(s_\alpha-\zeta): \alpha\in H^0\}\subseteq L\setminus (\textbf{V}_\theta+\textbf{U}_\theta).$$ Therefore $H_{\zeta}^0 \cap (\textbf{V}_\theta+\textbf{U}_\theta)=\varnothing,$ this shows that $\theta \notin \overline{H_{\zeta}^0}$ which is a contradiction.
	
{\bf (ii)} Let $\zeta$ be a \emph{rough weighted $\mathcal{I}_\tau$-cluster point} of $s=\{s_{\alpha}\}_{\alpha\in D}.$ Then for any  $\tau$-neighbourhood ${\textbf{U}_\theta}$ of $\theta$ we get $$A=\{\alpha\in D: t_\alpha(s_\alpha-\zeta)\in L\setminus (\textbf{V}_\theta+\textbf{U}_\theta)\}\notin \mathcal{F(\mathcal{I})}.$$ 
If $H \in \mathcal{F(\mathcal{I})}$ then $H \not\subseteq A.$ Then there exists $\alpha_0\in H$ such that $\alpha_0 \in D\setminus A$ (since $\alpha_0 \notin A$). Therefore $t_{\alpha_0}(s_{\alpha_0}-\zeta)\in H_{\zeta} \cap (\textbf{V}_\theta+\textbf{U}_\theta).$ Since $L$ is Hausdorff and $H_{\zeta} \cap (\textbf{V}_\theta+\textbf{U}_\theta)\neq \varnothing$ for all $\textbf{U}_\theta$   so $ H_{\zeta} \cap \overline{\textbf{V}_\theta} \neq \varnothing.$
Hence our assertion follows.
\end{proof}

We now derive a succinct characterisation  of closed set in terms of the set $\mathcal{R}_{t}^{\textbf{V}_\theta}-{\mathcal{I}}(\Gamma_{s}^{\tau}).$


\begin{theorem}  Let $(L, \tau)$ be a completely separable \emph{locally solid Riesz space}, $\textbf{U}_\theta$ be a closed $\tau$-neighbourhood of $\theta,$ $\mathcal{I}$ be a non-trivial ideal of a directed set $D$,  and there exist a pairwise disjoint sequence of sets $\{B_i\}_{i\in \mathbb{N}}$ such that $B_i \subset D,$ 	$B_i \notin \mathcal{I}$ for all $i\in \mathbb{N}.$   Then there is a net $s=\{s_{\alpha}\}_{\alpha\in D}$ in  $L$ and a weighted net $t=\{t_{\alpha}\}_{\alpha\in D}$ such that  $$\textbf{U}_\theta=\displaystyle{\bigcap_{\textbf{V}_\theta \in \mathcal{N}_{sol}}} \left(\mathcal{R}_{t}^{\textbf{V}_\theta}-{\mathcal{I}}(\Gamma_{s}^{\tau}) \right).$$    
	\end{theorem}

\begin{proof} Since $\textbf{U}_\theta$ is separable so there exists a countable dense set $A=\{a_1,a_2, \ldots\}$ of $ {\textbf{U}_\theta}$ such that $\overline{A}=\textbf{U}_\theta.$  Let  $s=\{s_{\alpha}\}_{\alpha\in D}$ be a net in the following manner:
	\begin{equation*}
	s_{\alpha}= \begin{cases} a_i  & \mbox {if} ~\alpha \in B_i, \\
	a_1 &  \mbox {if} ~  \alpha\in D\setminus \displaystyle{\bigcup_{i=1}^\infty} ~ B_i.
	\end{cases}
	\end{equation*}		
	We choose a  $\mathcal{I}$-bounded \emph{weighted net} $t=\{t_{\alpha}\}_{\alpha\in D}$ such that  $|t_\alpha|\geq \delta >1$ for all $\alpha\in D.$ From the  $\mathcal{I}$-boundedness of $t=\{t_{\alpha}\}_{\alpha\in D}$  we get a large positive real number $M$ such that 	 $$B=\{ \alpha \in D : |t_{\alpha}| \leq M \} \in \mathcal{F}(\mathcal{I}).$$	
	Let $u\in \textbf{U}_\theta.$ 	 Since $u+\frac{1}{M}\textbf{V}_\theta$ is a $\tau$-neighbourhood of $u$ for all $\textbf{V}_\theta \in \mathcal{N}_{sol}.$ Then
	$$A\cap (u+\frac{1}{M}\textbf{V}_\theta)\neq \varnothing ~\implies ~a_i\in u+\frac{1}{M}\textbf{V}_\theta ~\mbox{for some}~i\in \mathbb{N} $$ $$\implies ~M(a_i-u)\in \textbf{V}_\theta ~\mbox{for some}~i\in \mathbb{N}$$
	$$\implies ~t_{\alpha}(s_{\alpha}-u)\in \textbf{V}_\theta ~\mbox{where}~\alpha\in B\cap B_i~ \mbox{for some}~i\in \mathbb{N}$$
	$$\implies B\cap B_i \subseteq \{\alpha \in D: ~t_{\alpha}(s_{\alpha}-u)\in \textbf{V}_\theta +\textbf{W}_\theta\} \not\in \mathcal{I},$$
	where $\textbf{W}_\theta$ is any $\tau$-neighbourhood of $\theta$ (since $B\cap B_i \notin \mathcal{I},$ otherwise $B_i=(\mathbb{N}\setminus B)\cap (B\cap B_i)\in \mathcal{I},$ a contradiction). Therefore 	
	\begin{equation}\label{Closed-1} u \in \mathcal{R}_{t}^{\textbf{V}_\theta}-{\mathcal{I}}(\Gamma_{s}^{\tau})~\implies~ \textbf{U}_\theta \subseteq  \mathcal{R}_{t}^{\textbf{V}_\theta}-{\mathcal{I}}(\Gamma_{s}^{\tau}),
	~\mbox{where}~	\textbf{V}_\theta ~\mbox{is an arbitrary element of}~\mathcal{N}_{sol}.
	\end{equation}	
	Then,	
	\begin{equation}\label{Closed-2} \textbf{U}_\theta \subseteq  \displaystyle{\bigcap_{\textbf{V}_\theta \in \mathcal{N}_{sol}}} \left(\mathcal{R}_{t}^{\textbf{V}_\theta}-{\mathcal{I}}(\Gamma_{s}^{\tau})\right).
		\end{equation}
Without using the finite intersection property or compactness, we derive the condition $$ \displaystyle{\bigcap_{\textbf{V}_\theta \in \mathcal{N}_{sol}}} \left(\mathcal{R}_{t}^{\textbf{V}_\theta}-{\mathcal{I}}(\Gamma_{s}^{\tau})\right)\neq \varnothing$$ from equations \ref{Closed-1} and \ref{Closed-2}.\\


	Next, we assume that $$\zeta \in \displaystyle{\bigcap_{\textbf{V}_\theta \in \mathcal{N}_{sol}}} \left(\mathcal{R}_{t}^{\textbf{V}_\theta}-{\mathcal{I}}(\Gamma_{s}^{\tau}) \right).$$
	Also let $$C= \displaystyle{\bigcap_{\textbf{V}_\theta \in \mathcal{N}_{sol}}} \{\alpha\in D: t_\alpha(s_\alpha-\zeta)\in \textbf{V}_\theta+\textbf{V}_\theta\}.$$ If $C=\varnothing,$ there exist two $\tau$-neighbourhoods $\textbf{V}^{'}_\theta, \textbf{V}^{''}_\theta \in \mathcal{N}_{sol}$ such that 
	$$\{\alpha\in D: t_\alpha(s_\alpha-\zeta)\in \textbf{V}^{'}_\theta+\textbf{V}^{'}_\theta\} \cap \{\alpha\in D: t_\alpha(s_\alpha-\zeta)\in \textbf{V}^{''}_\theta+\textbf{V}^{''}_\theta\} =\varnothing.$$
	Since $2\textbf{V}^{'}_\theta \cap 2\textbf{V}^{''}_\theta$ is a $\tau$-neighbourhood of $\theta.$ Then there exists an element $\textbf{V}^{'''}_\theta \in \mathcal{N}_{sol}$ such that 
	$$\textbf{V}^{'''}_\theta + \textbf{V}^{'''}_\theta \subseteq 2\textbf{V}^{'}_\theta \cap 2\textbf{V}^{''}_\theta.$$ 
	Therefore $\{\alpha\in D: t_\alpha(s_\alpha-\zeta)\in \textbf{V}^{'''}_\theta+\textbf{V}^{'''}_\theta\} =\varnothing,$ which is a contradiction since  $$ \zeta \in \mathcal{R}_{t}^{\textbf{V}^{'''}_\theta}-{\mathcal{I}}(\Gamma_{s}^{\tau}) ~\mbox{and}~  \{\alpha\in D: t_\alpha(s_\alpha-\zeta)\in \textbf{V}^{'''}_\theta+\textbf{V}^{'''}_\theta\} \notin \mathcal{I}.$$	
	Let $\alpha_0\in C.$ Then $t_{\alpha_0}(s_{\alpha_0}-\zeta)\in \textbf{V}_\theta +\textbf{V}_\theta$ for all $\textbf{V}_\theta \in \mathcal{N}_{sol}.$ 	
	$$\implies ~\delta(a_i-\zeta)\in \textbf{V}_\theta+\textbf{V}_\theta~\mbox{for some}~i\in \mathbb{N},~\mbox{and}~~\mbox{for all}~ \textbf{V}_\theta \in \mathcal{N}_{sol}$$
	$$\implies ~a_i\in \zeta +\frac{1}{\delta}(\textbf{V}_\theta+\textbf{V}_\theta) \subseteq  \zeta + \textbf{V}_\theta+\textbf{V}_\theta ~\mbox{for some}~i\in \mathbb{N},~\mbox{and}~~\mbox{for all}~ \textbf{V}_\theta \in \mathcal{N}_{sol}.$$
	Therefore $a_i \in (\zeta +\textbf{V}_\theta+\textbf{V}_\theta) \cap \textbf{U}_\theta \neq \varnothing$  for all $\textbf{V}_\theta \in \mathcal{N}_{sol}.$ Since $\textbf{U}_\theta$ is a closed set so $\zeta \in \textbf{U}_\theta.$
	\end{proof}


In the year (2004/2005), \v{C}in\v{c}ura et al. \cite{CSST} give an important result related to the $\mathcal{I}$-cluster point set, which is stated below. 

\begin{theorem} \cite[Lemma 3.1]{CSST} Let $\mathcal{I}$ be an admissible ideal on $\mathbb{N},$ $\Gamma_{x}(\mathcal{I})$ be a set of all \emph{$\mathcal{I}$-cluster points of a sequence} $x=\{x_n\}_{n\in \mathbb{N}}$  and $K$ be a compact subset of a metric space $X.$ Then for every sequence $x=\{x_n\}_{n\in \mathbb{N}}$ in $X$ with $\{n\in \mathbb{N}: x_n\in K\}  \notin \mathcal{I}$ we have $$K \cap \Gamma_{x}(\mathcal{I})\neq \varnothing.$$  \end{theorem}
 
  While reformulating the above results based on the set ~`\emph{rough weighted $\mathcal{I}_\tau$-cluster points of a net}' in \emph{locally solid Riesz spaces}, the object of
 $$ K \cap \left(\mathcal{R}_{t}^{\textbf{V}_\theta}-{\mathcal{I}}(\Gamma_{s}^{\tau})\right)\neq \varnothing$$ is violated. We exemplify this assertion below.


\begin{example} \label{Example 4} Consider the \emph{locally solid Riesz space}, directed set $D,$ $G$ and ideal $\mathcal{I}$ as in Example \ref{Example 1} (i).    Also let the compact set $K=\{x\in \mathbb{R}^2 : ||x||\leq 1 \}=\{(x,y)\in \mathbb{R}^2: x^2+y^2\leq 1\}=\textbf{V}_\theta,$ net $s=\{s_{\alpha}\}_{\alpha\in D}$ and the weighted net $t=\{t_{\alpha}\}_{\alpha\in D}$ are in the following manner:
	\begin{equation*}
		s_{\textbf{U}_\theta(r)}= \begin{cases} (r,0)  & \mbox {if} ~\textbf{U}_\theta (r)\in G, \\
			(\frac{1}{r},\frac{1}{r}) & \mbox {if} ~  \textbf{U}_\theta (r)\in D\setminus G,
		\end{cases}
		~\mbox{and}~
		t_{\textbf{U}_\theta(r)}=\begin{cases} -(r+1) & \mbox {if}~  r\in G,\\
			2r+1 & \mbox {if}~  r\in D\setminus G.\\
		\end{cases}
	\end{equation*}
 Therefore $K \cap \{\mathcal{R}_{t}^{\textbf{V}_\theta}-{\mathcal{I}}(\Gamma_{s}^{\tau})\}= \varnothing$ and $\{\alpha \in D: s_\alpha \in K\}\not\in \mathcal{I}.$
	\end{example}


\begin{theorem} Let $\mathcal{I}$ be a non-trivial ideal of a directed set $D,$ $(L,\tau)$ be a \emph{locally solid Riesz space}, $\textbf{V}_\theta \in \mathcal{N}_{sol},$ $K$ be a compact subset of $L$ and the weighted net $t=\{t_{\alpha}\}_{\alpha\in D}$ be $\mathcal{I}$-bounded. Then for every net $s=\{s_{\alpha}\}_{\alpha\in D}$ in $L$ with $$\{\alpha \in D: s_\alpha\in K\}  \notin \mathcal{I}~\mbox{we have}~ K \cap \{\mathcal{R}^{\textbf{V}_\theta}_{t}-{\mathcal{I}}(\Gamma_{s}^{\tau})\}\neq \varnothing.$$	
\end{theorem}

\begin{proof} We assume by way of contradiction, let $K \cap \{\mathcal{R}^{\textbf{V}_\theta}_{t}-{\mathcal{I}}(\Gamma_{s}^{\tau})\}= \varnothing.$  Let $k\in K.$ Then $k\notin \mathcal{R}^{\textbf{V}_\theta}_{t}-{\mathcal{I}}(\Gamma_{s}^{\tau}),$ there exists a $\textbf{U}_\theta \in \mathcal{N}_{sol}$ such that
	 $$A_k=\{\alpha\in D: t_\alpha(s_\alpha-k)\in \textbf{V}_\theta+\textbf{U}_\theta\}\in \mathcal{I}.$$ From the \emph{$\mathcal{I}$-boundedness} of the weighted net $t=\{t_{\alpha}\}_{\alpha\in D}$ there exists $M>0$ such that $$B=\{\alpha\in D: |t_\alpha|\leq M\}\in \mathcal{F}(\mathcal{I}).$$  Also
	 $$K \subseteq  \displaystyle{\bigcup_{k\in K}}~\left(k+\frac{\textbf{V}_\theta}{M}\right).$$	 
Since $K$ is a compact set, then there exists a finite number of elements $k_1,k_2,\ldots,k_n$ in $K$ such that
$$K \subseteq  \displaystyle{\bigcup_{i=1}^n}~\left(k_i+\frac{\textbf{V}_\theta}{M}\right).$$
Let $C=\{\alpha\in D: s_\alpha \in K\}\notin \mathcal{I}.$ Then $B\cap C\neq \varnothing$ (otherwise $C\subseteq D\setminus B\in \mathcal{I},$ a contradiction). If $\alpha \in  B\cap C$ then
$$s_\alpha \in K \subseteq  \displaystyle{\bigcup_{i=1}^n}~\left(k_i+\frac{\textbf{V}_\theta}{M}\right)~\implies ~ s_\alpha \in k_i+\frac{\textbf{V}_\theta}{M}~\text{for some}~i\in \{1,2,3,\ldots,n\},$$
$$\implies ~ M(s_\alpha -k_i)\in \textbf{V}_\theta~\text{for some}~i\in \{1,2,3,\ldots,n\},$$
$$\implies ~ t_\alpha (s_\alpha -k_i)\in \textbf{V}_\theta~\text{for some}~i\in \{1,2,3,\ldots,n\},$$
	since $|t_\alpha(s_\alpha-k_i)|\leq |M(s_\alpha-k_i)|$ for all $i\in \{1,2,3,\ldots,n\},$ $\alpha \in B\cap C$ and $\textbf{V}_\theta \in \mathcal{N}_{sol}$. Therefore	$$B\cap C\subseteq \{\alpha\in D: t_\alpha(s_\alpha-k_i) \in \textbf{V}_\theta\} ~\text{for some}~i\in \{1,2,3,\ldots,n\},$$ 
$$\implies ~ B\cap C \subseteq \displaystyle{\bigcup_{i=1}^{n}}~\{\alpha\in D: t_\alpha(s_\alpha-k_i) \in \textbf{V}_\theta\}  \subseteq \displaystyle{\bigcup_{i=1}^{n}}~\{\alpha\in D: t_\alpha(s_\alpha-k_i) \in \textbf{V}_\theta + \textbf{U}_\theta\}=  \displaystyle{\bigcup_{i=1}^{n}}~ A_{k_i}\in \mathcal{I}.$$
Finally, $C\subseteq (D\setminus B)\cup (B\cap C)\in \mathcal{I},$ a contradiction. Therefore, this is the result.
\end{proof}


\begin{theorem}\label{Th 2.5}  Let $\mathcal{I}$ be a non-trivial ideal of a directed set $D,$ $t=\{t_{\alpha}\}_{\alpha\in D}$ be a \emph{weighted net}, $s=\{s_{\alpha}\}_{\alpha\in D}$ be a net in a \emph{locally solid Riesz space} $L,$ $\textbf{V}_\theta \in \mathcal{N}_{sol}$ and $\textbf{V}_\theta$ be a compact $\tau$-neighbourhood of $\theta.$ If $\mathcal{R}^{\textbf{V}_\theta}_{t}-{\mathcal{I}}(\Gamma_{s}^{\tau})\neq \varnothing$ then there exists a compact set $K$ of $L$ such that $$\{\alpha \in D: s_\alpha\in K\}  \notin \mathcal{I}.$$
\end{theorem}

\begin{proof} Let  $\zeta \in \mathcal{R}^{\textbf{V}_\theta}_{t}-{\mathcal{I}}(\Gamma_{s}^{\tau}).$ Then for every $\textbf{U}_\theta \in \mathcal{N}_{sol}$ and $\textbf{U}_\theta \subseteq  \textbf{V}_\theta$ such that
	$$\{\alpha\in D: t_\alpha(s_\alpha-\zeta)\in \textbf{V}_\theta+\textbf{U}_\theta\}\notin \mathcal{I}.$$
	$$\implies~ A=\{\alpha\in D: t_\alpha(s_\alpha-\zeta)\in \textbf{V}_\theta+\textbf{V}_\theta\}\notin \mathcal{I}.$$
	  Since $t=\{t_{\alpha}\}_{\alpha\in D}$ is a weighted net, there exists a positive real number $\delta$ such that  $|t_{\alpha}|\geq \delta$ for all  $\alpha\in D.$ Let $\alpha \in A.$ Then
	  $|\delta (s_\alpha-\zeta)| \leq |t_\alpha(s_\alpha-\zeta)|$ for all $\alpha \in A.$ Also, $\textbf{V}_\theta +\textbf{V}_\theta \in \mathcal{N}_{sol}$ and $t_\alpha(s_\alpha-\zeta)\in \textbf{V}_\theta+\textbf{V}_\theta$ implies $ \delta (s_\alpha-\zeta)\in \textbf{V}_\theta+\textbf{V}_\theta.$	  
	  Therefore 	
	$$A\subseteq \{\alpha\in D: \delta (s_\alpha-\zeta)\in \textbf{V}_\theta+\textbf{V}_\theta\}=
	 \{\alpha\in D:  s_\alpha\in K\}\notin \mathcal{I}$$ where $K=\zeta+\frac{1}{\delta}(\textbf{V}_\theta+\textbf{V}_\theta)$ and $K$ is a compact set in $L.$\end{proof}

The following example shows that if $\textbf{V}_\theta$ is not a compact $\tau$-neighbourhood of $\theta$ in the above theorem, then there does not exists any compact set $K$ of $L$ such that $\{\alpha \in D: s_\alpha\in K\}  \notin \mathcal{I}.$

\begin{example} Let us consider the infinite dimensional norm linear space $l^2$ space with the norm
$$||x||=\sqrt{\left(\displaystyle{\sum_{i=1}^{\infty}}~|x_i|^2 \right)}~\mbox{where}~x=\{x_{n}\}_{n\in \mathbb{N}}\in l^2$$ with coordinatewise ordering. Since  $l^2$ is a \emph{locally solid Riesz space}. The family $\mathcal{N}_{sol}$ of all $\textbf{U}_\theta(r)$ defined as $\textbf{U}_\theta(r)=\{x\in l^2 :||x||<r\}$ where $r>0,$ constitutes a base at $\theta=(0,0,\ldots).$ The set $D=\mathcal{N}_{sol}$ is a directed set with the binary relation $\geq$ defined by $\textbf{U}_\theta\geq \textbf{V}_\theta$ if $\textbf{U}_\theta\supseteq \textbf{V}_\theta$ for any $\textbf{U}_\theta,\textbf{V}_\theta\in D.$ Consider the following decomposition of $\mathbb{R}^+,$  

$$B_0=(0,1)\cup \{\displaystyle{\bigcup_{n\in O}}~[n,n+1]\}~\mbox{where} ~O=\{1,3,5,\ldots\},$$

$$B_\ell=\displaystyle{\bigcup_{n\in O}}~(2^{\ell} n, 2^{\ell} n+1)~\mbox{where} ~ \ell\in \mathbb{N}=\{1,2,3,\ldots\},$$

	$$\implies~\mathbb{R}^+=\displaystyle{\bigcup_{\ell \in \mathbb{N}_0}}~B_\ell~\mbox{where}~~B_i \cap B_j=\varnothing ~\mbox{for all}~ i,j\in \mathbb{N}_0=\{0,1,2,3,\ldots\} ~\mbox{and}~i\neq j.$$

Let $\mathcal{I}=\left\{\{\textbf{U}_\theta(r): r\in \Delta (\subset \mathbb{R})\}: \Delta ~\mbox{can intersect finitely many elements of}~ B_i,~i\in \mathbb{N}_0\right\}.$ We choose the net $s=\{s_{\alpha}\}_{\alpha\in D}$ and the weighted net $t=\{t_{\alpha}\}_{\alpha\in D}$ are respectively:
\begin{equation*}		s_{{\textbf{U}_\theta}(r)}=\begin{cases} \theta &  \mbox {if} ~ r\in B_0, \\
	e_n & \mbox {if} ~             ~~~~~ r \in B_n,
	\end{cases}
	~\mbox{and}~
	t_{U(r)}= \begin{cases} r+1 & \mbox {if} ~ r\in B_0, \\
	1 & \mbox {if} ~             ~~~~~ r \in B_n,
	\end{cases}
\end{equation*}
where $e_{n}=(\delta^n_{j})$ has $n^{\text{th}}$ term $1$ and all other terms are $0.$ Let $\textbf{V}_\theta=\{x\in l^2 :||x||\leq 1\}.$ For any $\tau$-neighbourhood $\textbf{U}_\theta$ of $\theta$ the set $ \{\alpha \in D: t_\alpha(s_\alpha-\theta)\in \textbf{V}_\theta+\textbf{U}_\theta\}=D \notin \mathcal{I}.$
Therefore $\theta \in \mathcal{R}^{\textbf{V}_\theta}_{t}-{\mathcal{I}}(\Gamma_{s}^{\tau})\neq \varnothing.$ \\

Next, let $K$ be any subset of $L.$ Then two possibility are their
\begin{description}
	\item[(i)] Let $K\cap (\{e_i:i\in \mathbb{N}\}\cup \{\theta\})=\{e_{m_1},e_{m_2},e_{m_3},\ldots,e_{m_k}\}$ be a finite set. Then 	
    $$\{\alpha \in D: s_\alpha\in K\}=\{ \textbf{U}_\theta(r): r\in \displaystyle{\bigcup_{i=1}^k}~B_{m_i}\}  \in \mathcal{I}.$$
	
	\item[(ii)] Let $K\cap (\{e_i:i\in \mathbb{N}\}\cup \{\theta\})=\{e_{m_1},e_{m_2},e_{m_3},\ldots\}$ be an infinite set. Then the set 
	$$\{\alpha \in D: s_\alpha\in K\}=\{ \textbf{U}_\theta(r): r\in \displaystyle{\bigcup_{i=1}^\infty}~B_{m_i}\}  \notin \mathcal{I}.$$
\end{description}

Therefore the only possibility the set $\{\alpha \in D: s_\alpha\in K\}  \notin \mathcal{I}$ if $K$ contains infinitely many $e_i$ where $i\in \mathbb{N}.$    But this is not possible for any compact set $K$ of $l^2.$


 \end{example}


In general if $\textbf{V}_\theta$ be \emph{$\tau$-bounded} (or \emph{$\tau$-unbounded}) then the set $\mathcal{R}_{t}^{\textbf{V}_\theta}-{\mathcal{I}}(\Gamma_{s}^{\tau})$  may be  \emph{$\tau$-unbounded}. To prove this assertion, we consider the following examples:


\begin{example} {\bf Case (i):} Let $\textbf{V}_\theta$ be \emph{$\tau$-bounded}.  To show that the set $\mathcal{R}_{t}^{\textbf{V}_\theta}-{\mathcal{I}}(\Gamma_{s}^{\tau})$ is  \emph{$\tau$-unbounded}. This example based on Example \ref{Example 2} (i). Replacing the net $s=\{s_{\alpha}\}_{\alpha\in D}$ by the net $s^*=\{s^*_{\alpha}\}_{\alpha \in D}$ where
	\begin{equation*}
	s^*_{{\textbf{U}_\theta}(r)}=\begin{cases} (1,1)   & \mbox {if}  ~ r\in B_0, \\
	(\ell+\frac{1}{r^2},1)  & \mbox {if} ~r\in B_{\ell}, ~\mbox{where}~\ell\in \mathbb{N},
	\end{cases}
	\end{equation*}
	and keeping all other conditions unaltered in Example \ref{Example 2}, consequently we get $$\mathbb{N}\times \{1\} \subseteq \mathcal{R}_{t}^{\textbf{V}_\theta}-\mathcal{I}(\Gamma_{s^*}^{\tau}).$$
	
	{\bf Case (ii):}  Let $\textbf{V}_\theta$ be \emph{$\tau$-unbounded}. Consider the \emph{locally solid Riesz space}, directed set $D,$ a subset $G$ of $D$ and ideal $\mathcal{I}$ as in Example \ref{Example 1} (i). Let  the net $s=\{s_{\alpha}\}_{\alpha\in D}$ and the weighted net $t=\{t_{\alpha}\}_{\alpha\in D}$ be in the following manner:
	\begin{equation*}
	s_{\textbf{U}_\theta(r)}\in \begin{cases} (r,0) & \mbox {if} ~\textbf{U}_\theta (r)\in G, \\
	(4,4)+ (4 \cos r, 4\sin r) & \mbox {if}~ \textbf{U}_\theta (r)\in D\setminus G,
	\end{cases}
	~\mbox{and}~ ~ ~
	t_{\textbf{U}_\theta(r)}=\begin{cases} -(r+1) & \mbox {if}~  r\in G,\\
	r+1 & \mbox {if}~ r\in D\setminus G.\\
	\end{cases}
	\end{equation*}
	Let $\textbf{V}_\theta=\{(x,y)\in \mathbb{R}^2 :-1<y<1\}.$ Then $\mathcal{R}_{t}^{\textbf{V}_\theta}-{\mathcal{I}}(\Gamma_{s}^{\tau})=\{(x,y)\in \mathbb{R}^2 :y=0\}.$  This example also indicates the net $s=\{s_{\alpha}\}_{\alpha\in D}$ is not \emph{rough weighted  $\mathcal{I}_\tau$-convergent} to any point of $\mathbb{R}^2$ w.r.t the roughness of degree $\textbf{V}_\theta$  yet the set $\mathcal{R}_{t}^{\textbf{V}_\theta}-{\mathcal{I}}(\Gamma_{s}^{\tau})$ is uncountable and \emph{$\tau$-unbounded}. \end{example}

We give some elementary properties concerning $\tau$-boundedness of the set $\mathcal{R}^{\textbf{V}_\theta}_{t}-{\mathcal{I}}(\Gamma_{s}^{\tau}).$


\begin{prop} Let $(L,\tau)$ be a locally solid Riesz space, $s=\{s_{\alpha}\}_{\alpha\in D}$ be  weighted  $\mathcal{I}_\tau$-bounded,  $t=\{t_{\alpha}\}_{\alpha\in D}$ be a weighted net, 	$T=\displaystyle{\inf_{\alpha\in D}} |t_\alpha|$ and ${\textbf{V}}_\theta\in \mathcal{N}_{sol}.$ Then
	\begin{description}
		\item[(i)] $\mathcal{R}^{\textbf{V}_\theta}_{t}-{\mathcal{I}}(\Gamma_{s}^{\tau})\subseteq \frac{2+\mu}{T} ~ {\textbf{V}_\theta},$ where $\mu$ is a positive real number,\\ 
	\item[(ii)]$\mathcal{R}^{\textbf{V}_\theta}_{t}-{\mathcal{I}}(\Gamma_{s}^{\tau})$ is $\tau$-bounded if  $\textbf{V}_\theta$ is $\tau$-bounded.	
	\end{description}	
		 \end{prop}

\begin{proof} Since $s=\{s_{\alpha}\}_{\alpha\in D}$ is \emph{weighted  $\mathcal{I}_\tau$-bounded}, then for every $\tau$-neighbourhood ${\textbf{V}_\theta}$ of $\theta$ there exists some $\lambda>0$ such that  $A=\{\alpha\in D: \lambda t_\alpha s_\alpha \notin {\textbf{V}_\theta}\}\in \mathcal{I}.$ If $\zeta\in \mathcal{R}^{\textbf{V}_\theta}_{t}-{\mathcal{I}}(\Gamma_{s}^{\tau}),$ then $B=\{\alpha\in D: t_\alpha (s_\alpha-\zeta) \in {\textbf{V}_\theta}+{\textbf{V}_\theta}\}\notin \mathcal{I}.$ Also, $(D\setminus A)\cap B \neq \varnothing$ (otherwise $B\subseteq A\in \mathcal{I},$ a contradiction). If $\alpha \in (D\setminus A)\cap B,$ then $$t_\alpha \zeta=t_\alpha (\zeta-s_\alpha)+t_\alpha s_\alpha \in (2+\frac{1}{\lambda}) ~ \textbf{V}_\theta.$$
 Since ${\textbf{V}}_\theta \in \mathcal{N}_{sol}$ then  $T \zeta\in  (2+\frac{1}{\lambda}) ~ \textbf{V}_\theta.$ Hence $\mathcal{R}^{\textbf{V}_\theta}_{t}-{\mathcal{I}}(\Gamma_{s}^{\tau}) \subseteq  \frac{2+\mu}{T} ~{\textbf{V}_\theta}$ where $\mu=\frac{1}{\lambda}.$\\

	For the next part  $\textbf{V}_\theta$ is $\tau$-bounded so $\frac{2+\mu}{T} ~{\textbf{V}_\theta}$ also.		 
	\end{proof}
	
As an {\bf application} of the sets 	
	$ \text{W}\mathcal{I}_{\tau}-\text{LIM}^{\textbf{V}_\theta}_t s$ and  $\mathcal{R}^{\textbf{V}_\theta}_{t}-{\mathcal{I}}(\Gamma_{s}^{\tau})$  we prove the following result which characterise the ideal.


\begin{theorem}\label{Th 2.7}  Let $(L,\tau)$ be a \emph{Hausdorff locally solid Riesz space}, $\mathcal{I}$ be a non-trivial ideal of a directed set $D.$  Then $\mathcal{I}$ is a maximal ideal on $D$ if and only if 
	$ \text{W}\mathcal{I}_{\tau}-\text{LIM}^{\textbf{V}_\theta}_t s=\mathcal{R}^{\textbf{V}_\theta}_{t}-{\mathcal{I}}(\Gamma_{s}^{\tau})$ for any net $s=\{s_{\alpha}\}_{\alpha\in D}$ in  $L,$  with some roughness degree $\textbf{V}_\theta.$
\end{theorem}

\begin{proof} Let $\mathcal{I}$ be a maximal ideal on $D.$ From Definition \ref{Def 2.4} and \ref{Def 2.5} we get 	$$\text{W}\mathcal{I}_{\tau}-\text{LIM}^{\textbf{V}_\theta}_t s\subseteq \mathcal{R}^{\textbf{V}_\theta}_{t}-{\mathcal{I}}(\Gamma_{s}^{\tau})~\mbox{for any roughness degree}~ \textbf{V}_\theta \in \mathcal{N}_{sol}.$$
	
	Next let, $\zeta \in \mathcal{R}^{\textbf{V}_\theta}_{t}-{\mathcal{I}}(\Gamma_{s}^{\tau}).$ Then $A=\{\alpha \in D: t_\alpha (s_\alpha- \zeta)\in \textbf{V}_\theta+\textbf{U}_\theta \}\notin \mathcal{I},$ for any $\tau$-neighbourhood $\textbf{U}_\theta$ of $\theta.$ By the property of maximal ideal we get $D\setminus A\in \mathcal{I}.$ This implies 
	$$D\setminus A=\{\alpha \in D: t_\alpha (s_\alpha- \zeta)\notin \textbf{V}_\theta+\textbf{U}_\theta \}\in \mathcal{I}.$$
	So $\zeta \in \text{W}\mathcal{I}_{\tau}-\text{LIM}^{\textbf{V}_\theta}_t s.$ Therefore 	$ \text{W}\mathcal{I}_{\tau}-\text{LIM}^{\textbf{V}_\theta}_t s=\mathcal{R}^{\textbf{V}_\theta}_{t}-{\mathcal{I}}(\Gamma_{s}^{\tau}).$
	
	Conversely, let $\text{W}\mathcal{I}_{\tau}-\text{LIM}^{\textbf{V}_\theta}_t s=\mathcal{R}^{\textbf{V}_\theta}_{t}-{\mathcal{I}}(\Gamma_{s}^{\tau})$ for any net $s=\{s_{\alpha}\}_{\alpha\in D}$ in  $L,$  with some roughness degree $\textbf{V}_\theta.$  With out any loss of generality we assume that $\mathcal{I}$ is not a maximal ideal on $D.$ Then there exists a subset $A$ of a directed set $D$ such that $A$ and $D\setminus A$ are both not in $\mathcal{I}.$ Let $x\neq \theta.$ Then there exists a $\tau$-neighbourhood $\textbf{V}^*_\theta\in \mathcal{N}_{sol}$ such that $\theta \in 2 \textbf{V}^*_\theta$ but $Tx \notin 2 \textbf{V}^*_\theta,$ where $T=\displaystyle{\inf_{\alpha\in D}} |t_\alpha|.$	Let  $s=\{s_{\alpha}\}_{\alpha\in D}$ be a net in the following manner:
	\begin{equation*}
	s_{\alpha}= \begin{cases} \theta  & \mbox {if} ~\alpha \in A, \\
	x &  \mbox {if} ~  \alpha\in D\setminus A.
	\end{cases}
	\end{equation*}
	
	Since $\textbf{V}^*_\theta\in \mathcal{N}_{sol}$ and $T\leq |t_\alpha|$ for all $\alpha \in D$ and $Tx \notin 2 \textbf{V}^*_\theta$ implies $t_\alpha x \notin 2\textbf{V}^*_\theta.$ Then	
	 $$A=\{\alpha \in D: t_\alpha (s_\alpha- \theta)\in \textbf{V}^*_\theta+\textbf{V}^*_\theta \}$$
	 $$=\{\alpha \in D: t_\alpha (s_\alpha- \theta)\in \textbf{V}^*_\theta+\textbf{U}_\theta \}, ~\mbox{where}~ \textbf{U}_\theta \subseteq \textbf{V}^*_\theta ~\mbox{and}~ \textbf{U}_\theta \in \mathcal{N}_{sol},$$ 	
	$$\implies D\setminus A=\{\alpha \in D: t_\alpha (s_\alpha- x)\in \textbf{V}^*_\theta+\textbf{V}^*_\theta \}$$
	$$=\{\alpha \in D: t_\alpha (s_\alpha- x)\in \textbf{V}^*_\theta+\textbf{U}_\theta \}, ~\mbox{where}~\textbf{U}_\theta \subseteq \textbf{V}^*_\theta~\mbox{and}~\textbf{U}_\theta \in \mathcal{N}_{sol}.$$ 
	Since $A$ and $D\setminus A$ not in $\mathcal{I}.$ Then $x, \theta \in \mathcal{R}^{\textbf{V}^*_\theta}_{t}-{\mathcal{I}}(\Gamma_{s}^{\tau})=\text{W}\mathcal{I}_{\tau}-\text{LIM}^{\textbf{V}^*_\theta}_t s.$ Hence  $A$ and $D\setminus A$ are both lies in $\mathcal{F}(\mathcal{I}),$ a contraction since $\varnothing \notin \mathcal{F}(\mathcal{I}).$ So $\mathcal{I}$ is an maximal ideal on $D.$	\end{proof}

\begin{Corollary}  Let $\mathcal{I}$ be a non-trival ideal on the directed set $D,$ $t=\{t_{\alpha}\}_{\alpha\in D}$ be a \emph{weighted net}, $s=\{s_{\alpha}\}_{\alpha\in D}$ be a net in a \emph{locally solid Riesz space} $L.$ If $\left(\mathcal{R}_{t}^{\textbf{V}_\theta}-{\mathcal{I}}(\Gamma_{s}^{\tau}) \right) \setminus \left( \text{W}\mathcal{I}_{\tau}-\text{LIM}^{\textbf{V}_\theta}_t s \right) \neq \varnothing$ then $\mathcal{I}$ is not a maximal ideal of $D.$
\end{Corollary}

\begin{proof} Let $\zeta \in \left(\mathcal{R}_{t}^{\textbf{V}_\theta}-{\mathcal{I}}(\Gamma_{s}^{\tau}) \right) \setminus \left( \text{W}\mathcal{I}_{\tau}-\text{LIM}^{\textbf{V}_\theta}_t s \right).$ Then $\{\alpha \in D: t_\alpha (s_\alpha- \zeta)\in \textbf{V}_\theta+\textbf{U}_\theta \}\notin \mathcal{I},$ for any $\tau$-neighbourhood $\textbf{U}_\theta$ of $\theta.$ 
	Since $\zeta \notin \text{W}\mathcal{I}_{\tau}-\text{LIM}^{\textbf{V}_\theta}_t s.$ Therefore $\{\alpha \in D: t_\alpha (s_\alpha- \zeta)\notin \textbf{V}_\theta+\textbf{U}_\theta \}\notin \mathcal{I}.$ Hence the ideal $\mathcal{I}$ is not a maximal ideal on $D.$
\end{proof}

\section{{Declarations}}

{\bf Ethical Approval} : Not applicable.

\vspace{0.2cm}

{\bf Authors' Contributions} : All the authors have equally contributed to this article.

\vspace{0.2cm}

{\bf Funding} : No funding.

\vspace{0.2cm}

{\bf Data Availibility} : This manuscript has no associated data.

\vspace{0.2cm}

{\bf Competing interests} :  The authors declare that they have no known competing financial interests or personal relationships that could have appeared to influence the work reported in this paper. 



\begin{thebibliography}{99}

\large




\bibitem{AP} H. Albayrak, S. Pehlivan, Statistical convergence and statistical continuity on locally solid Riesz spaces, Topology Appl. {\bf 159} (2012) 1887--1893.




\bibitem{A1} {S. Aytar}, The rough limit set and the core of a real sequence, Numer. Funct. Anal. Optim. {\bf 29} (3-4) (2008)  283--290.

\bibitem{A2} {S. Aytar}, Rough statistical convergence, Numer. Funct. Anal. Optim. {\bf 29} (3--4) (2008) 291--303.

\bibitem{A3} {S. Aytar}, Rough statistical cluster points, Filomat {\bf 31} (16) (2017) 5295--5304.


\bibitem{Aziz-Ghosal} T. Aziz,  S. Ghosal, $f$-rough Cauchy sequences, Quaest. Math. \textbf{48} (7) (2025) 1047--1063.


\bibitem{BL} M. Balcerzak, P. Leonetti, On the relationship between ideal cluster points and ideal limit points, Topology Appl. {\bf 252} (2019) 178--190.


\bibitem{CSST} {J. \v{C}in\v{c}ura, T. \v{S}al\'{a}t,  M. Sleziak,  V. Toma}, Sets of statistical cluster points and $\mathcal{I}$-cluster points, Real Anal. Exchange {\bf 30} (2004/2005)  565--580.



\bibitem{Das-1} P. Das, Some further results on ideal convergence in topological spaces, Topology Appl. {\bf 159} (2012) 2621--2626.


\bibitem{Das-Savas2014(Taiwan)} P. Das, E. Sava\c{s},  On $\mathcal{I}$-statistical pre-Cauchy sequences, Taiwan. J. Math. {\bf 18} (1) (2014) 115--126.



\bibitem{DS(filomat)}  {P. Das, E. Savas}, On $\mathcal{I}$-convergence of nets in locally solid Riesz spaces, Filomat {\bf 27} (1) (2013)  89--94.


\bibitem{Di} G. Di Maio,  Li. D. R. Ko$\check{\mbox{c}}$inac,  Statistical convergence in topology,  Topology Appl. {\bf 156} (2008) 28--45.



\bibitem{Dundar2014} E. D$\ddot{\mbox{u}}$ndar, C. {\c{C}}akan, Rough $\mathcal{I}$-convergence, Demonstr. Math. \textbf{47} (3) 638--651 (2014).


\bibitem{Dundar2016} E. D$\ddot{\mbox{u}}$ndar, On rough $\mathcal{I}_2$-convergence of double sequences, Numer. Funct. Anal. Optim. {\bf 37} (4) (2016) 480--491.




\bibitem{Fridy}  {J. A. Fridy},  On statistical convergence, Analysis {\bf 5} (4) (1985) 301--313.

\bibitem{fridy1993statistical}  J. A. Fridy, Statistical limit points,  Proc. Amer. Math. Soc. {\bf 118} (4) (1993) {1187--1192}.

\bibitem{GB} {S. Ghosal, M. Banerjee}, Rough weighted statistical convergence on locally solid Riesz spaces, Positivity {\bf 25} (2021) 1789--1804.

\bibitem{G3} S. Ghosal, M. Banerjee, A. Ghosh, Weighted modulus $S_{\theta}$-convergence of order $\alpha$ in probability, Arab. J.  Math. Sci.  \textbf{23} (2017) 242--257.


\bibitem{GM1} {S. Ghosal, S. Mandal}, The degree of roughness, Topology Appl. {\bf 307} (2022)  107944.

\bibitem{GM2} {S. Ghosal, S. Mandal}, Rough weighted $\mathcal{I}$-$\alpha\beta$-statistical convergence in locally solid Riesz spaces, J. Math. Anal. Appl. {\bf 506} (2022) 125681.






\bibitem{Kos1} P. Kostyrko, M. Ma\v{c}aj, T. \v{S}al\'{a}t, M. Sleziak, $\mathcal{I}$-convergence and extremal limit points,  Math. Slovaca {\bf 55} (4) (2005) 443--464.



\bibitem{Kos2} P. Kostyrko, T. \v{S}al\'{a}t, O. Strauch, $\mathcal{I}$-convergence, Real Anal. Exchange {\bf 26} (2) (2000/2001) 669--686.

\bibitem{Lahiri} B. K. Lahiri, P. Das, $\mathcal{I}$ and $\mathcal{I}^*$-convergence of nets, Real Anal. Exchange {\bf 33} (2) (2007/2008) 431--442.

\bibitem{Leonetti} P. Leonetti, Rough Families, Cluster Points, and Cores, J. Convex Anal. {\bf 32} (4) (2025) 1083--1090.


\bibitem{LR1} M. C. List$\acute{\mbox{a}}$n-Garc\'ia, F. Rambla-Barreno, A characterization of uniform rotundity in every direction in terms of rough convergence, Numer. Funct. Anal. Optim. {\bf 32} (11) 1166--1174 (2011).

\bibitem{Garcia} M. C. List$\acute{\mbox{a}}$n-Garc\'ia, F. Rambla-Barreno, Rough convergence and Chebyshev centers in Banach spaces, Numer. Funct. Anal. Optim. {\bf 35} (4) (2014) 432--442.




\bibitem{Mandal-Ghosal} S. Mandal, S. Ghosal,  When rough $\mathcal{I}$-convergent sequences fail to be rough $\mathcal{I}$-Cauchy, Hacet. J. Math. Stat. {\bf Doi}: 10.15672/hujms.1837190 (2026).


\bibitem{pal} S. K. Pal, D. Chandra, S. Dutta,  Rough ideal convergence, Hacet. J. Math. Stat. {\bf 42} (6) (2013) 633--640.


\bibitem{Phu1} {H. X. Phu}, Rough convergence in normed linear spaces, Numer. Funct. Anal. Optim. {\bf 22} (1-2) (2001) 199--222.

\bibitem{Phu2}  {H. X. Phu}, Rough convergence in infinite dimensional normed space, Numer. Funct. Anal. Optimiz. {\bf 24} (3-4) (2003) 285--301.

\bibitem{mursaleen} S. A. Rahaman, M. Mursaleen, On rough deferred statistical convergence of difference sequences in $\mathcal{L}$-fuzzy normed spaces, J. Math. Anal. Appl. {\bf 530} (2) (2024) 127684.

\bibitem{R}  {F. Riesz}, Sur la decomposition des operations fonctionelles lin\'{e}aires, Alti del Congr. Internaz. del Mat. Bologna 1928, 3, Zanichelli (1930) 143--148.

\bibitem{Salat} {T. \v{S}al\'{a}t},  On statistical convergent sequences of real numbers, Math. Slovaca {\bf 30} (2) (1980) 139--150.


\bibitem{SGSG(Taiwan)} A. Sahiner, M. Gurdal, S. Saltan, H. Gunawan, Ideal convergence in $2$-normed spaces, Taiwan. J. Math. {\bf 11} (5) (2007) 1477--1484.

\bibitem{SP} C. Sencimen, S. Pehlivan, Statistical order convergence in Riesz spaces, Math. Slovaca {\bf 62} (2012) 257--270.




\bibitem{Steinhaus} H. Steinhaus, Sur la convergence ordinaire et la convergence asymptotique. Colloq. Math. {\bf 2} (1951) 73--74.

\end{thebibliography}
\end{document}